\documentclass{article}

\usepackage{PRIMEarxiv}

\usepackage[utf8]{inputenc} % allow utf-8 input
\usepackage[T1]{fontenc}    % use 8-bit T1 fonts

\usepackage{url}            % simple URL typesetting
\usepackage{booktabs}       % professional-quality tables
\usepackage{amsmath,amssymb,amsfonts}
\usepackage{amsthm}
\usepackage{cite}
\usepackage{nicefrac}       % compact symbols for 1/2, etc.
\usepackage{microtype}      % microtypography
\usepackage{lipsum}
\usepackage{fancyhdr}       % header
\usepackage{graphicx}       % graphics
\graphicspath{{media/}}     % organize your images and other figures under media/ folder

\title{Adaptive Regulation with Global KL Guarantees}

\author{
  Iasson Karafyllis \\
   Department of Mathematics \\
  National Technical University of Athens \\
   Zografou Campus, 15780, Athens, Greece\\
  \texttt{iasonkar@central.ntua.gr} \\
   \And
 Alexandros Aslanidis\\
   Department of Mathematics \\
  National Technical University of Athens \\
   Zografou Campus, 15780, Athens, Greece\\
  \texttt{alex\textunderscore aslanidis@mail.ntua.gr} \\
   \And
  Miroslav Krstic \\
   Department of Mechanical and Aerospace Engineering \\
  University of California\\
 San Diego, La Jolla, CA 92093 USA\\
  \texttt{krstic@ucsd.edu} \\
}

\newtheorem{thm}{Theorem}
\newtheorem{lem}{Lemma}
\newtheorem{defin}{Definition}
\newtheorem{assumption}{Assumption}

\begin{document}
\maketitle

\begin{abstract}
In the absence of persistency of excitation (PE), referring to adaptive control systems as ``asymptotically stable’’ typically indicates insufficient understanding of stability concepts. While the state is indeed regulated to zero and the parameter estimate has some limit, namely, the overall state converges to some equilibrium, the equilibrium reached is not unique (and not even necessarily stable) but is dependent on the initial condition. The equilibrium set in the absence of PE is not uniformly attractive (from an open set containing the equilibrium set); hence, asymptotic stability does not hold and KL estimates are unavailable for the full state. In this paper we pursue adaptive control design with KL guarantees on the regulated state, something that is possible but previously unachieved with smooth, time-invariant and non-hybrid adaptive controllers. This property is referred to as Uniform Global Asymptotic Output Stability, where the regulated state is thought of as a system output. We provide designs for (i) systems with a matched uncertainty and (ii) systems in the parametric strict feedback form. To guarantee KL estimates in the absence of PE, our designs employ time-invariant nonlinear damping terms, which depend both on the state and the parameter estimate. With an example, we illustrate the theory.
\end{abstract}

\section{Introduction}
The problem of adaptive control arises naturally as it is rare (or even impossible) to know exactly the system parameters. Thus, it is of great interest to regulate a system in the presence of unknown parameters. The literature on adaptive control is huge, spanning over sixty years. It is well known that adaptive control does not guarantee uniform attractivity and/or KL-estimates for the solution of the closed-loop system (see \cite{b1}) when the dynamics are time-invariant and when no Persistence of Excitation (PE) assumption holds. Indeed, as shown in \cite{bKr} the equilibrium set in the absence of PE is not uniformly attractive (from an open set containing the equilibrium set); hence, asymptotic stability does not hold and KL estimates are unavailable for the full state. In the literature one can find a number of different approaches for the problem of constructing adaptive feedback laws that can provide some performance guarantees for the closed-loop system:
\begin{enumerate}
\item[i)]
Using the funnel control scheme (see for instance \cite{b17}, \cite{b30}, \cite{b31}, \cite{b32}, \cite{b33}) researchers have constructed time-varying adaptive controllers that guarantee a prescribed performance in a region of the state space. Similar to the funnel control scheme, the prescribed performance scheme (see \cite{b16}, \cite{b34}, \cite{b35}, \cite{b36}, \cite{b15}, \cite{b37}) employs an error transformation of the system and guarantees prescribed performance locally.
\item[ii)]
For tracking control problems researchers have managed to establish Global Uniform Attractivity under the assumption of PE (see \cite{b13}, \cite{b14}, \cite{b18}, \cite{b19} and also \cite{b40}). Using the Immersion and Invariance approach researchers in \cite{b27} managed to construct a strict Lyapunov function for a class of nonlinearly parameterized systems under an ISS small-gain condition and a PE assumption, and guarantee Global Uniform Asymptotic Tracking. The Immersion and Invariance approach is also used in  \cite{b25}, \cite{b28}.
\item[iii)]
In \cite{b20}, \cite{b21} an unbounded-in-time, time-varying adaptive control scheme is proposed by using a time-depended state transformation in order to guarantee Uniform Attractivity properties in the absence of PE.
\item[iv)]
In \cite{b22},\cite{b23} an event-triggered adaptive control scheme is proposed in order to guarantee specific KL-estimates for the solutions of the closed-loop system.
\item[v)] In \cite{b54f} an adaptive scheme is proposed that guarantees regulation of the state in finite time uniformly, assuming that a bound for the parameter vector is known, see also \cite{b52f}, \cite{b53f} for finite time regulation using state-feedback gains. An unbounded time-scale transformation is used in \cite{b51f} in order to regulate the state to the origin in finite time and similarly, time-varying gains are used in \cite{b50f} to guarantee global prescribed time convergence of the state to zero.   

\end{enumerate}

It should also be noticed that the recent work \cite{b5} showed that minor modifications of the adaptive controller design procedures introduced in \cite{b1} can also guarantee Semi-Global Uniform Output Stability for nonlinear time-invariant systems under a smooth, time-invariant, adaptive feedback law.

In this paper we approach the problem of design with KL guarantees of the regulated state by applying adaptive backstepping in \cite{b1} using nonlinear damping terms that depend not only on the state, but also on the estimate of the parameter. We apply the control design on systems satisfying the matching condition (see \cite{b1}) and derive a KL-estimate which guarantees the Uniform Global Asymptotic Output Stability (UGAOS) for a class of output functions. This design follows the idea in \cite{b5} where semi-global results were provided. Then we extend the design to systems in the parametric strict feedback form, where we derive a KL-estimate that quantifies the convergence speed to a residual ball of assignable radius for the regulated state and yields UGAOS for a similar class of outputs as in the previous case.

The article is organized as follows. In Section~\ref{sec:A} we present the controllers and prove the desired properties in the case of systems satisfying the matching condition and in Section~\ref{sec:B} we consider the case of systems in the parametric strict feedback form. In Section~\ref{sec:Example} we give an illustrative example of the design and we compare the performance of the proposed adaptive controller with standard adaptive feedback laws proposed in \cite{b1}.
\section*{Notation and notions}
In what follows we use the notation $(x)^{+}$ for the positive part of the real number $x \in \mathbb{R}$, i.e. $(x)^{+}= \max{(x,0)}$. For a vector $x \in \mathbb{R}^{n}$, $\vert x \vert$ denotes the Euclidean norm of $x$ and $x'$ the transpose of $x$.

Let $f:\mathbb{R}^{n} \rightarrow \mathbb{R}^{n}$ be a locally Lipschitz vector field with $f(0)=0$ and $h:\mathbb{R}^{n} \rightarrow \mathbb{R}^{k}$ be a continuous mapping with $h(0)=0$. Consider the dynamical system
\begin{equation}
\dot{x}=f(x), x \in \mathbb{R}^{n} \label{ds}
\end{equation}
with output
\begin{equation}
Y=h(x) \label{dso}
\end{equation}
We assume that the dynamical system \eqref{ds} is forward complete, i.e., for every $x_{0} \in \mathbb{R}^{n}$ the unique solution $x(t)=\phi(t,x_{0})$ of the initial value problem \eqref{ds} with initial condition $x(0)=x_{0}$ exists for all $t \geq 0$.

We use the notation $Y(t,x_{0})=h(\phi(t,x_{0}))$ for all $t \geq 0$, $x_{0} \in \mathbb{R}^{n}$ and $B_{R}=\left\lbrace x \in \mathbb{R}^{n}: \vert x \vert <R \right\rbrace$ for all $R>0$. The following properties are standard in the analysis of output stability: see for instance \cite{b4}.

\begin{defin}[\textbf{Output Stability Notions}] \label{def:0}
We say that the system \eqref{ds}, \eqref{dso} is
\begin{enumerate}
\item[i)]Lagrange output stable if for every $R>0$ the set $\left\lbrace \vert Y(t,x_{0})\vert: x_{0} \in B_{R}, t \geq 0 \right\rbrace$ is bounded.
\item[ii)] Lyapunov output stable if for every $\epsilon>0$ there exists $\delta(\epsilon)>0$ such that for all $x_{0} \in B_{\delta(\epsilon)}$, it holds that $\vert y(t,x_{0})\vert \leq \epsilon$ for all $t \geq 0$.
\item[iii)]Globally Asymptotically Output Stable (GAOS) if system \eqref{ds}, \eqref{dso} is Lagrange and Lyapunov output stable and $\lim_{t \rightarrow +\infty}(Y(t,x_{0}))$ for all $x_{0} \in \mathbb{R}^{n}$.
\item[iv)] Uniformly Globally Asymptotically Output Stable (UGAOS) if system \eqref{ds}, \eqref{dso} is Lagrange and Lyapunov output stable and for every $\epsilon,R>0$ there exists $T(\epsilon,R)>0$ such that for all $x_{0} \in B_{R}$, it holds that $\vert Y(t,x_{0})\vert \leq \epsilon$ for all $t \geq T(\epsilon,R)$.
\end{enumerate}

\end{defin} 

When the word ''output'' is omitted in all the above stability notions then we understand that the output coincides with the (regulated) state of the plant. The following lemma provides a characterization of the UGAOS notion of Definition 1. Its proof is given in \cite{b3}.

\begin{lem}[\textbf{UGAOS through a KL estimate}] System \eqref{ds}, \eqref{dso} is UGAOS if and only if there exists a function $\sigma \in KL$ such that the following estimate holds for all $x_{0} \in \mathbb{R}^{n}$ and $t \geq 0$:
\begin{equation}
\vert Y(t,x_{0})\vert \leq \sigma (\vert x_{0} \vert,t)
\end{equation}
\end{lem}

\section{Systems satisfying the matching condition} \label{sec:A}
Consider the following finite-dimensional control system

\begin{equation}
\dot{x}=f(x)+g(x)u+g(x)(\varphi(x))'\theta 
\label{erl}
\end{equation}
where $f,g: \mathbb{R}^n \rightarrow \mathbb{R}^n , \varphi :\mathbb{R}^n \rightarrow \mathbb{R}^p$ are locally Lipschitz vector fields with $f(0)=0$ and  $\varphi (0)=0$ , $x \in \mathbb{R}^n$, $u \in \mathbb{R}$ is the control input and $\theta \in \mathbb{R}^p$ is the constant vector of unknown parameters. The control system \eqref{erl} satisfies the "matching condition" (following the terminology in \cite{b1}) since the uncertain parameters are in the span of the control.
The following assumption is used for the control system \eqref{erl}.

\begin{assumption}\label{as H}
There exist functions $P \in C^2 \left(  \mathbb{R}^n ; \mathbb{R}_{+} \right) , Q \in C^0 \left(  \mathbb{R}^n ; \mathbb{R}_{+} \right) $ being positive definite and radially unbounded and a locally Lipschitz function $k \in C^0 \left(  \mathbb{R}^n ; \mathbb{R}_{+} \right)$ with $k(0)=0$ such that 
\begin{equation}
\nabla P(x)f(x)+\nabla P(x)g(x)k(x)\leq -Q(x) , \text{for all } x \in \mathbb{R}^n \label{hh}
\end{equation}
\end{assumption}
 When Assumption \ref{as H}  holds then it is possible to regulate the state to zero by applying a standard adaptive control scheme. To see this, consider the adaptive controller 

\begin{align}
u &= k(x)- \left( \varphi (x) \right)' \hat{\theta} \label{eq:sac0} \\
\dfrac{d\hat{\theta}}{dt} &= \Gamma \nabla P(x) g(x) \varphi (x) \label{eq:sac}
\end{align}
where $\Gamma \in \mathbb{R}^{p \times p}$ is a positive definite matrix. The closed-loop system \eqref{erl}, \eqref{eq:sac0}, \eqref{eq:sac} with output 

\begin{equation}
Y=x \label{outp}
\end{equation}
has the following stability properties:

\begin{enumerate}
\item[i)] the set $\left\lbrace \left( 0,\hat{\theta} \right) : \hat{\theta} \in \mathbb{R}^p \right\rbrace $ (i.e., the set where $x=0$) is the set of equilibrium points for the closed-loop system \eqref{erl},\eqref{eq:sac}
\item[ii)] the equilibrium point $\left( x,\hat{\theta} \right) =\left( 0, \theta \right) $ is Lyapunov stable
\item[iii)]for every initial condition, the solution $\left(  x(t),\hat{\theta(t) }  \right) $ of the closed-loop system \eqref{erl},\eqref{eq:sac} is bounded and the closed-loop system is Lagrange stable
\item[iv)] the closed-loop system \eqref{erl},\eqref{eq:sac},\eqref{outp} is GAOS and every solution of the closed-loop system \eqref{erl},\eqref{eq:sac} satisfies $ \lim_{t\rightarrow +\infty} \left( x(t) \right) =0 $

\end{enumerate}
 The proof of the above facts follows the analysis given in \cite{b1} using the Lyapunov function $V(x,\hat{\theta}) := P(x)+\dfrac{1}{2} \left( \hat{\theta}-\theta \right)' \Gamma^{-1} \left( \hat{\theta}-\theta \right)$, which (by virtue \eqref{erl},\eqref{hh}) satisfies the following differential inequality
 
\begin{equation}
\dfrac{d}{dt} V(x,\hat{\theta}) \leq -Q(x) \label{difen}
\end{equation}
However, the convergence rate of $x \in \mathbb{R}^{n}$ to zero is not necessarily uniform with respect to the initial condition, even when considering initial states in a bounded neighbourhood of the origin.
In order to be able to guarantee a KL-estimate we need an additional assumption.
\smallskip
\begin{assumption} \label{as A}
There exists a locally Lipschitz, positive function $\mu \in C^{0}\left( \mathbb{R}^{n} ; (0,+\infty) \right)$ such that the following inequality holds:
\begin{equation}
\vert \varphi(x) \vert^{2} \leq \mu (x)Q(x) , \text{ for all }x \in \mathbb{R}^{n} \label{aa}
\end{equation}
\end{assumption}

Assumption \ref{as A} is automatically guaranteed if there exists a constant $q>0$ such that the inequality $Q(z) \geq q \vert z \vert^{2}$ holds for all $z \in \mathbb{R}^{n}$ in a neighbourhood of $0 \in \mathbb{R}^{n}$. Both Assumption \ref{as H} and Assumption \ref{as A} are automatically guaranteed if the feedback law $u=k(x)$ achieves global asymptotic stabilization and local exponential stabilization of $0 \in \mathbb{R}^{n}$ for the control system $\dot{x}=f(x)+g(x)u$.

When Assumption \ref{as H} and Assumption \ref{as A} hold then it is possible to modify the adaptive control scheme \eqref{eq:sac0}, \eqref{eq:sac} so that we can obtain two different KL-estimates: a "perfect" KL-bound when $\vert \theta \vert \leq \sqrt{r}$ and a "practical" KL-bound when $\vert \theta \vert > \sqrt{r}$, where $r \geq 0$ is a parameter of the controller. More specifically, the adaptive control scheme \eqref{eq:sac0}, \eqref{eq:sac} is modified in the following way:

\begin{align}
u =& k(x)-(\varphi(x))'\hat{\theta} \left( \dfrac{\delta}{2} \vert \varphi(x) \vert^{2}+ \mu (x) \left( r+\vert \hat{\theta} \vert^{2} \right)\right)\left( \nabla P(x)g(x) \right) \label{eq:mac0} \\
\dfrac{d\hat{\theta}}{dt} =& \Gamma \nabla P(x) g(x) \varphi(x) \label{eq:mac}  
\end{align}
where $\Gamma \in \mathbb{R}^{p \times p}$ is a positive definite matrix and $\delta>0$, $r \geq 0$ are constants. The term $- \left( \dfrac{\delta}{2} \vert \varphi(x) \vert^{2}+ \mu (x) \left( r+\vert \hat{\theta} \vert^{2} \right)\right)\left( \nabla P(x)g(x) \right)$ is the difference between \eqref{eq:sac0}, \eqref{eq:sac} and \eqref{eq:mac0}, \eqref{eq:mac} and it is a nonlinear damping term. Nonlinear damping terms have been proposed in the literature (see Section 3.4.3 in \cite{b1}) and it was shown in \cite{b5} that nonlinear damping terms can guarantee semi-global stabilization. Here the nonlinear damping term has a crucial difference with the damping term proposed in \cite{b1} (and in \cite{b5}): the nonlinear damping term does not depend on $x$ only (as in \cite{b1} and \cite{b5}) but on both $x$ and $\hat{\theta}$. That is the reason that here we achieve global and uniform stability properties. The controller parameter $r \geq 0$ plays a special role in our design: we achieve UGAOS for the output $Y=x$ when $\vert \theta \vert \leq \sqrt{r}$. 

It should be noticed that all properties (i), (ii), (iii) and (iv) hold for the closed-loop system \eqref{erl} with \eqref{eq:mac0}, \eqref{eq:mac}. The proof of this fact follows again the analysis given in \cite{b1} using the Lyapunov function $V(x,\hat{\theta}) := P(x)+ \dfrac{1}{2} \left( \hat{\theta}-\theta \right)' \Gamma^{-1} \left( \hat{\theta}-\theta \right) $ , which satisfies the differential inequality \eqref{difen}. However, the closed-loop system \eqref{erl} with \eqref{eq:mac0}, \eqref{eq:mac} also satisfies a global KL-estimate given by the following theorem. 

\begin{thm} \label{thm:1}
Suppose that Assumptions \ref{as H} and \ref{as A} hold. Let $\rho \in K_{\infty}$ be a function for which the following inequality holds:
\begin{equation}
\rho \left(P(x) \right) \leq Q(x) \text{ for all } x \in \mathbb{R}^{n} \label{eq:t1h} 
\end{equation} 
Then for every $ \lambda \in \left(0,1 \right)$, there exists $\sigma \in KL$ such that the following global estimate holds for every solution of the closed-loop system \eqref{erl} with \eqref{eq:mac}:
\begin{equation}
\left( P(x(t))-\alpha \right)^{+} \leq \sigma \left( \left( P(x(0))-\alpha\right)^{+} ,t \right) \text{ for all t } \geq 0 \label{eq:gkl1}  
\end{equation}
where
\begin{equation}
\alpha := \rho^{-1}\left( (1-\lambda)^{-1} \delta^{-1}\left( \vert \theta \vert^{2} -r \right)^{+}  \right) \label{eq:defa}
\end{equation}
\end{thm}

\textbf{Remarks:} \begin{enumerate}
\item[(a)] Since $P \in C^{2}(\mathbb{R}^{n};\mathbb{R}_{+})$ and $Q \in C^{0}(\mathbb{R}^{n};\mathbb{R}_{+})$ are positive definite and radially unbounded functions, there exists $\rho \in K_{\infty}$ such that inequality \eqref{eq:t1h} holds. Therefore, the existence of $\rho \in K_{\infty}$ for which \eqref{eq:t1h} holds is not an assumption.
\item[(b)] Since $P \in C^{2}(\mathbb{R}^{n};\mathbb{R}_{+})$ is a positive definite and radially unbounded function, there exist functions $b_{1},b_{2} \in K_{\infty}$ such that
\begin{equation}
b_{1}(\vert x \vert) \leq P(x) \leq  b_{2}(\vert x \vert) , \text{for all x} \in \mathbb{R}^{n} \label{eq:pdru}
\end{equation}
The proof of Theorem \ref{thm:1} and inequality \eqref{eq:gkl1} shows that the following global KL-estimate holds for all $t \geq 0$:
\begin{equation}
b_{1}(\vert x(t) \vert) \leq \tilde{\sigma} \left( \left( b_{2}(\vert x(0)\vert)-\alpha \right)^{+},t  \right)+\alpha \label{eq:gklr1} 
\end{equation}
for appropriate $\tilde{\sigma} \in KL$. It should be noticed that the functions $b_{1},b_{2} \in K_{\infty}$ and $\tilde{\sigma} \in KL$ are independent of $\theta$. Only $\alpha := \rho^{-1}\left( (1-\lambda)^{-1} \delta^{-1}\left( \vert \theta \vert^{2} -r \right)^{+}  \right)$ depends on $\theta$. Furthermore, notice that by increasing $\delta>0$ or $r>0$ (which are parameters of the controller), we can decrease the value of $\alpha$.
\item[(c)] The KL-estimate \eqref{eq:gkl1} guarantees UGAOS with any output $Y=h(x,\theta)$ with $h(0,\theta)=0$ for all $\theta \in \mathbb{R}^{p}$ for which there exists a function $b\in K_{\infty}$ with $\vert h(x,\theta)\vert \leq b \left( \left( P(x)-\alpha \right)^{+} \right)$ for all $(x,\theta) \in \mathbb{R}^{n} \times \mathbb{R}^{p}$. 
\item[(d)] Estimate \eqref{eq:gklr1} guarantees UGAOS for the output $Y=x$ when $\vert \theta \vert \leq \sqrt{r}$. On the other hand, estimate \eqref{eq:gklr1} is a practical KL-bound when $\vert \theta \vert > \sqrt{r}$: it guarantees that $x(t)$ will approach a neighbourhood of radius $b_{1}^{-1}(\alpha)$ in a uniform (with to the initial conditions) rate.
\end{enumerate}
Inequality \eqref{eq:gklr1} shows the process of ''simultaneous learning and doing''. Getting to the vicinity of the ball $\vert x \vert \leq b^{-1}_{1}(\alpha)$ can be fast and occurs in a uniform (with respect to the initial conditions) rate. However, if the system is unstable and highly uncertain (large $\vert \theta \vert$), this vicinity will be large, and the remaining task (getting to $x=0$) will occur at an unquantified rate. This happens because when we approach the set $x=0$ there is no excitation to learn, and to then speed up the remainder of the convergence as a result of learning. 

\section{The parametric strict feedback form} \label{sec:B}

In this section we study the parametric strict feedback form given by the following equations for $i=1,\ldots,n$ :

\begin{align} \label{eq:csthl26}
\dot{x}_{i}=f_{i}(x_{1},\ldots,x_{i})+g_{i}(x_{1},\ldots,x_{i})x_{i+1}+\sum_{j=1}^{p}\varphi_{i,j}(x_{1},\ldots,x_{i})\theta_{j} 
\end{align}
where $f_{i},g_{i}:\mathbb{R}^{i}\rightarrow \mathbb{R}$, $\varphi_{i,j}:\mathbb{R}^{i}\rightarrow \mathbb{R}$, $i=1,\ldots,n$, $j=1,\ldots,p$ are $C^{\infty}$ functions with $f_{i}(0)=0$, $\varphi_{i,j}(0)=0$, $g_{i}(x_{1},\ldots,x_{i}) \neq 0$ for all $x \in \mathbb{R}^{n}$, $\theta \in \left( \theta_{1}, \theta_{2}, \ldots,\theta_{p} \right)' \in \mathbb{R}^{p}$ is a vector of unknown constant parameters and $x_{n+1}=u \in \mathbb{R}$ is the control input.

In order to study system \eqref{eq:csthl26} we need to introduce a useful notion.

\begin{defin}  \label{def:1}
A smooth parameterized family of diffeomorphisms of $\mathbb{R}^{n}$ is a $C^{\infty}$ map $T:\mathbb{R}^{n} \times \mathbb{R}^{p}\rightarrow \mathbb{R}^{n}$ with $T(0,\theta)=0$ for all $\theta \in \mathbb{R}^{p}$ having the property that there exists a $C^{\infty}$ map $T^{-1}:\mathbb{R}^{n} \times \mathbb{R}^{p}\rightarrow \mathbb{R}^{n}$ with $T\left( T^{-1}(x,\theta)\right)=T^{-1}\left( T(x,\theta)\right)=x $ for all $(x,\theta) \in \mathbb{R}^{n} \times \mathbb{R}^{p}$.
\end{defin}

The following technical results are also useful for our analysis.

\begin{lem} \label{lem1}
Let $h:\mathbb{R}^{n} \times \mathbb{R}^{p}\rightarrow \mathbb{R}$ be a $C^{\infty}$ map that satisfies $h(0,\theta)=0$ for all $\theta \in \mathbb{R}^{p}$. Then there exists a $C^{\infty}$ map $\rho :\mathbb{R}^{n} \times \mathbb{R}^{p}\rightarrow [1,+\infty)$ such that $\vert h(x,\theta)\vert \leq \rho ((x,\theta))\vert x \vert$ for all $(x,\theta) \in \mathbb{R}^{n} \times \mathbb{R}^{p}$.
\end{lem}

\begin{lem} \label{lem2}
Let $h:\mathbb{R}^{n} \times \mathbb{R}^{p}\rightarrow \mathbb{R}$ be a $C^{\infty}$ map that satisfies $h(0,\theta)=0$ for all $\theta \in \mathbb{R}^{p}$. Let $T:\mathbb{R}^{n} \times \mathbb{R}^{p}\rightarrow \mathbb{R}^{n}$ be a smooth parameterized family of diffeomorphisms of $\mathbb{R}^{n}$. Then there exists a $C^{\infty}$ map $\rho :\mathbb{R}^{n} \times \mathbb{R}^{p}\rightarrow [1,+\infty)$ such that $\vert h(x,\theta)\vert \leq \rho ((x,\theta))\vert T(x,\theta) \vert$ for all $(x,\theta) \in \mathbb{R}^{n} \times \mathbb{R}^{p}$.
\end{lem}

\begin{lem} \label{lem3}
Let $T:\mathbb{R}^{n} \times \mathbb{R}^{p}\rightarrow \mathbb{R}^{n}$ be a smooth parameterized family of diffeomorphisms of $\mathbb{R}^{n}$. Then for every bounded set $\Theta \subset \mathbb{R}^{p}$ there exist $K_{\infty}$ functions $\alpha, \beta : \mathbb{R}_{+} \rightarrow \mathbb{R}_{+}$ such that $\alpha(\vert x \vert) \leq \vert T(x,\theta) \vert \leq \beta(\vert x \vert)$ for all $(x,\theta) \in \mathbb{R}^{n} \times \Theta$. 
\end{lem}

The following lemma is the main backstepping lemma that allows us to “add integrators”. 

\begin{lem} \label{lem4}
Consider the system
\begin{equation}
\dot{x}=F(x)+b(x)y+\Phi(x)\theta \label{eq:cslemma}
\end{equation}
where $x \in \mathbb{R}^{n}$, $F,b:\mathbb{R}^{n} \rightarrow \mathbb{R}^{n}$ are $C^{\infty}$ vector fields with $F(0)=0$, $y \in \mathbb{R}$, $\Phi(x)$ is a $n \times p$ matrix with entries $C^{\infty}$ functions $\varphi_{i,j}:\mathbb{R}^{j} \rightarrow \mathbb{R}$ with $\varphi_{i,j}(0)=0$ for all $i=1,\ldots,n$ $j=1,\ldots,p$ and $\theta \in \mathbb{R}^{p}$ is the vector of unknown constant parameters. Suppose that there exist $C^{\infty}$ functions $k:\mathbb{R}^{n} \times \mathbb{R}^{p}\rightarrow \mathbb{R}$, $w:\mathbb{R}^{n} \times \mathbb{R}^{p}\rightarrow \mathbb{R}^{p}$ with $k(0,\hat{\theta})=0$ and $w(0,\hat{\theta})=0$ for all $\hat{\theta}\in \mathbb{R}^{p}$, a smooth parameterized family of diffeomorphisms of $\mathbb{R}^{n}$ $T:\mathbb{R}^{n} \times \mathbb{R}^{p}\rightarrow \mathbb{R}^{n}$, constants $r \geq 0$, $\alpha,\omega,\epsilon >0$ and $\gamma_{j}>0$, $j=1,\ldots,p$ such that the closed-loop system \eqref{eq:cslemma} with the adaptive feedback law:
\vspace{-0mm}
\begin{align}
\dfrac{d\hat{\theta}}{dt} &=w(x,\hat{\theta}) \label{eq:afllem0} \\
y &=k(x,\hat{\theta}) \label{eq:afllem}
\end{align}
satisfies the following inequalities for all $(x,\hat{\theta},\theta)\in \mathbb{R}^{n} \times \mathbb{R}^{p} \times \mathbb{R}^{p}$:

\begin{equation}
\dfrac{d}{dt}V(x,\hat{\theta}) \leq -\alpha \vert T(x,\hat{\theta}) \vert^{2} \label{eq:inV}
\end{equation}
\vspace{-0mm}
\begin{equation}
\dfrac{d}{dt}\left( \dfrac{1}{2}\vert T(x,\hat{\theta}) \vert^{2} \right) \leq -\omega \vert T(x,\hat{\theta}) \vert^{2}+ \epsilon \left( \vert \theta \vert^{2}-r \right)^{+} \label{eq:inrt} 
\end{equation}
where
\begin{equation}
V(x,\hat{\theta}):= \dfrac{1}{2}\vert T(x,\hat{\theta}) \vert^{2} + \sum_{j=1}^{p} \dfrac{1}{2\gamma_{j}}\left( \hat{\theta}_{j}-\theta_{j} \right)^{2} \label{eq:defV} 
\end{equation} 
Consider the control system \eqref{eq:cslemma} with 

\begin{equation}
\dot{y}=f(x,y)+g(x,y)u+ \sum_{j=1}^{p} \varphi_{n+1,j}(x,y)\theta_{j} \label{eq:2eqcslem}
\end{equation}
where $u \in \mathbb{R}$ the control input, $f,g:\mathbb{R}^{n+1} \rightarrow \mathbb{R}$, $\varphi_{n+1,j}:\mathbb{R}^{n+1} \rightarrow \mathbb{R}$ are $C^{\infty}$ functions with $f(0,0)=\varphi_{n+1,j}(0,0)=0$ for all $j=1,\ldots,p$ and $g(x,y) \neq 0$ for all $(x,y) \in \mathbb{R}^{n+1}$. Then, there exist $C^{\infty}$ functions $\tilde{k} : \mathbb{R}^{n+1} \times \mathbb{R}^{p}\rightarrow \mathbb{R}$, $\tilde{w} : \mathbb{R}^{n+1} \times \mathbb{R}^{p}\rightarrow \mathbb{R}^{p}$ with $\tilde{k}(0,0,\hat{\theta})=0$ and $\tilde{w}(0,0,\hat{\theta})=0$ for all $\hat{\theta}\in \mathbb{R}^{p}$, a smooth patameterized family of diffeomorphisms of $\mathbb{R}^{n+1}$ $\tilde{T} : \mathbb{R}^{n+1} \times \mathbb{R}^{p}\rightarrow \mathbb{R}^{n+1}$ such that the closed-loop system \eqref{eq:cslemma}, \eqref{eq:2eqcslem} with the adaptive feedback law   
\vspace{-0.5mm}
\begin{align}
\dfrac{d\hat{\theta}}{dt} &=\tilde{w}(x,y,\hat{\theta}) \label{eq:afl2lem0} \\
u &=\tilde{k}(x,y,\hat{\theta}) \label{eq:afl2lem}
\end{align}
satisfies the following inequalities for all $(x,y,\hat{\theta},\theta)\in \mathbb{R}^{n+1} \times \mathbb{R}^{p} \times \mathbb{R}^{p}$

\begin{equation}
\dfrac{d}{dt}\tilde{V}(x,y,\hat{\theta}) \leq -\alpha \vert \tilde{T}(x,y,\hat{\theta}) \vert^{2} \label{eq:inV2}
\end{equation}

\begin{equation}
\dfrac{d}{dt}\left( \dfrac{1}{2}\vert \tilde{T}(x,y,\hat{\theta}) \vert^{2} \right) \leq -\omega \vert \tilde{T}(x,y,\hat{\theta}) \vert^{2}+ 2\epsilon \left( \vert \theta \vert^{2}-r \right)^{+} \label{eq:inrt2} 
\end{equation}
where
\begin{equation}
\tilde{V}(x,y,\hat{\theta}):= \dfrac{1}{2}\vert \tilde{T}(x,y,\hat{\theta}) \vert^{2} + \sum_{j=1}^{p} \dfrac{1}{2\gamma_{j}}\left( \hat{\theta}_{j}-\theta_{j} \right)^{2} \label{eq:defV2} 
\end{equation}  

\end{lem}

\textbf{Remark on Lemma \ref{lem4}:}
Lemma \ref{lem4} is the analogue of Lemma 4.7 on page 134 in \cite{b1}. However, the proof and the control design of Lemma \ref{lem4} is different from the proof and the control design of Lemma 4.7 in \cite{b1}. Lemma \ref{lem4} exploits the idea of transient performance improvement described in Section 3.4.3 in \cite{b1}. More specifically, instead of using nonlinear damping terms (or state-dependent controller gains) that depend on $x$ only (as done in Section 3.4.3 in \cite{b1}), here we use nonlinear damping terms that depend both on $x$ and $\hat{\theta}$. Thus, we guarantee two differential inequalities (inequalities \eqref{eq:inV2} and \eqref{eq:inrt2}). Inequality \eqref{eq:inrt2} is the analogue of inequality (3.143) on page 118 in \cite{b1} with an essential difference: in inequality (3.143) of \cite{b1} the perturbation term depends on $\vert \hat{\theta}-\theta \vert$ while in \eqref{eq:inrt2} the perturbation term depends on $\theta$ only. As in the previous section, the two differential inequalities \eqref{eq:inV2}, \eqref{eq:inrt2} play a different role. Inequality \eqref{eq:inV2} is used in order to guarantee all properties that  adaptive control schemes usually guarantee (properties (i)-(iv) below). Inequality \eqref{eq:inrt2} provides a "perfect" KL-bound when $\vert \theta \vert \leq \sqrt{r}$ (i.e. UGAOS for the output $Y=x$) and a practical KL-bound when $\vert \theta \vert > \sqrt{r}$. As the reader can easily see these objectives are achieved by the simultaneous use of two Lyapunov-like functions: the function $\tilde{V}(x,\hat{\theta})$, defined by \eqref{eq:defV2}, which is the classical Lyapunov function used in adaptive control, and the function $\dfrac{1}{2} \vert \tilde{T}(x,\hat{\theta}) \vert^2$ which is a sort of an Input-to-Output Stability (IOS) function. This "dyad" of Lyapunov-like functions is an innovation of our analysis. Inequality \eqref{eq:inrt2} is the analogue of inequality \eqref{eq:gkl1} in Theorem \ref{thm:1}.    

\begin{thm} \label{thm:2}
Consider the control system \eqref{eq:csthl26}. Let arbitrary constants $r \geq 0$, $\alpha,\omega,\epsilon >0$ and $\gamma_{j}>0$, $j=1,\ldots,p$ be given. Then there exist $C^{\infty}$ functions $k:\mathbb{R}^{n} \times \mathbb{R}^{p} \rightarrow \mathbb{R}$, $w:\mathbb{R}^{n} \times \mathbb{R}^{p} \rightarrow \mathbb{R}^{p}$ with $k(0,\hat{\theta})=0$, $w(0,\hat{\theta})=0$ for all $\hat{\theta} \in \mathbb{R}^{p}$ 
and a smooth parameterized family of diffeomorphisms $T:\mathbb{R}^{n} \times \mathbb{R}^{p} \rightarrow \mathbb{R}^{n}$ such that the closed-loop system \eqref{eq:csthl26} with the adaptive feedback law 
\begin{align}
\dfrac{d\hat{\theta}}{dt} &=w(x,\hat{\theta}) \label{eq:th2ac0} \\
u &=k(x,\hat{\theta}) \label{eq:th2ac}
\end{align}
satisfies the following inequalities for all $(x,\hat{\theta},\theta)\in \mathbb{R}^{n} \times \mathbb{R}^{p} \times \mathbb{R}^{p}$:

\begin{equation}
\dfrac{d}{dt}V(x,\hat{\theta}) \leq -\alpha \vert T(x,\hat{\theta}) \vert^{2} \label{eq:th2inV}
\end{equation}

\begin{equation}
\dfrac{d}{dt}\left( \dfrac{1}{2}\vert T(x,\hat{\theta}) \vert^{2} \right) \leq -\omega \vert T(x,\hat{\theta}) \vert^{2}+ \epsilon \left( \vert \theta \vert^{2}-r \right)^{+} \label{eq:th2inrt} 
\end{equation}
where $V:\mathbb{R}^{n} \times \mathbb{R}^{p} \rightarrow \mathbb{R}$ is defined by \eqref{eq:defV}. Moreover, the following properties hold:

\begin{enumerate}
\item[(i)]the set $\left\lbrace (0,\hat{\theta}):\hat{\theta}\in \mathbb{R}^{p} \right\rbrace$ is the set of equilibrium points for the closed-loop system \eqref{eq:csthl26} with \eqref{eq:th2ac0}, \eqref{eq:th2ac},
\item[(ii)]the equilibrium point $(x,\hat{\theta})=(0,\theta)$ is Lyapunov stable for the closed-loop system \eqref{eq:csthl26} with \eqref{eq:th2ac0}, \eqref{eq:th2ac},
\item[(iii)]for every initial condition the solution $(x(t),\hat{\theta}(t))$ is bounded and system \eqref{eq:csthl26} with \eqref{eq:th2ac0}, \eqref{eq:th2ac} is Lagrange stable,
\item[(iv)]the closed-loop system \eqref{eq:csthl26} with \eqref{eq:th2ac0}, \eqref{eq:th2ac} and output $Y=x$ is GAOS and every solution satisfies $\lim_{t \rightarrow \infty}(x(t))=0$,
\item[(v)]the closed-loop system \eqref{eq:csthl26} with \eqref{eq:th2ac0}, \eqref{eq:th2ac} and output $Y=x$ is UGAOS when $\vert \theta \vert \leq \sqrt{r}$,
\item[(vi)]the closed-loop system \eqref{eq:csthl26} with \eqref{eq:th2ac0}, \eqref{eq:th2ac} is UGAOS for every output $Y=h(x,\hat{\theta},\theta)$ for which there exists a function $b \in K_{\infty}$ with  $\vert h(x,\hat{\theta},\theta) \vert \leq b\left( \left( \vert T(x,\hat{\theta}) \vert^{2}-\omega^{-1} \epsilon \left( \vert \theta \vert^{2}-r \right)^{+} \right)^{+} \right)$.
\end{enumerate}
\end{thm}

\textbf{Remark:}
All properties (i)-(iv) can be also guaranteed by the adaptive controllers proposed in \cite{b1}. However, here we also guarantee the uniform stability properties (v)-(vi). The proof of Theorem \ref{thm:2} shows that every solution of the closed-loop system satisfies the following estimate for all $t \geq 0$
\begin{align}
 \left( \vert T(x(t),\hat{\theta}(t)) \vert^{2}-\omega^{-1} \epsilon \left( \vert \theta \vert^{2}-r \right)^{+} \right)^{+} \leq \exp(-2\omega t) \left( \vert T(x(0),\hat{\theta}(0)) \vert^{2}-\omega^{-1} \epsilon \left( \vert \theta \vert^{2}-r \right)^{+} \right)^{+} 
\end{align}
The above estimate is an exponential estimate (a KL-estimate of exponential type) for $\left( \vert T(x(t),\hat{\theta}(t)) \vert^{2}-\omega^{-1} \epsilon \left( \vert \theta \vert^{2}-r \right)^{+} \right)^{+}$. Therefore, the above estimate shows again the process of ''simultaneous learning and doing''. Getting to the vicinity of the region $\vert T(x,\hat{\theta}) \vert \leq \sqrt{\omega^{-1} \epsilon \left( \vert \theta \vert^{2}-r \right)^{+}}$ can be fast and occurs in a uniform (with respect to the initial conditions) rate. However, as in the previous case of systems that satisfy a matching condition, if the system is unstable and highly uncertain (large $\vert \theta \vert$), this vicinity will be large, and the remaining task (getting to $x=0$) will occur at an un-quantified rate.

\section{Example} \label{sec:Example}
\vspace{-0mm}
Consider the simplified Moore-Greitzer model of jet engine surge dynamics,
\vspace{-0mm}
\begin{align}
\dot{x}_{1} &= \theta_{1}x_{1}^{2}+\theta_{2}x_{1}^{3}+x_{2} \label{eq:excs0} \\
\dot{x}_{2} &= u \label{eq:excs}
\end{align}
where $x=(x_{1},x_{2})' \in \mathbb{R}^{2}$ is the state, $u \in \mathbb{R}$ is the control input and $\theta=(\theta_{1},\theta_{2})' \in \mathbb{R}^{2}$ is the vector of constant unknown parameters. The nominal values for the parameters are $\theta=\left( -\dfrac{3}{2},-\dfrac{1}{2}\right)$ (see \cite{b1}) but since model \eqref{eq:excs0}, \eqref{eq:excs} is a simplified model obtained by Garlekin decomposition with only three modes, it is safer to consider $\theta=(\theta_{1},\theta_{2})' \in \mathbb{R}^{2}$ to be unknown.

Using the design methodology presented in \cite{b1}, we obtain the following adaptive feedback law for \eqref{eq:excs0}, \eqref{eq:excs}:
\vspace{-0mm}
\begin{align}
u =& -(Q^{-1}+\mu^{2})x_{1}-w_{1}x_{1}^{2}-w_{2}x_{1}^{3}   -(2\hat{\theta}_{1}x_{1}+3\hat{\theta}_{2}x_{1}^{2}+2\mu)(\hat{\theta}_{1}x_{1}^{2}+\hat{\theta}_{2}x_{1}^{3}+x_{2}) \label{exackkk00} \\
\dfrac{d\hat{\theta}_{1}}{dt} =& w_{1} = \gamma_{1}Qx_{1}^{2}(x_{2}+\hat{\theta}_{1}x_{1}^{2}+\hat{\theta}_{2}x_{1}^{3}+\mu x_{1}) (2\hat{\theta}_{1}x_{1}+3\hat{\theta}_{2}x_{1}^{2}+\mu)+\gamma_{1}x_{1}^{3} \label{exackkk0} \\
\dfrac{d\hat{\theta}_{2}}{dt} =& w_{2} = \gamma_{2}Qx_{1}^{3}(x_{2}+\hat{\theta}_{1}x_{1}^{2}+\hat{\theta}_{2}x_{1}^{3}+\mu x_{1})(2\hat{\theta}_{1}x_{1}+3\hat{\theta}_{2}x_{1}^{2}+\mu)+\gamma_{2}x_{1}^{4} \label{exackkk}
\end{align}
where $Q,\gamma_{1},\gamma_{2},\mu >0$ are constants. The adaptive feedback law \eqref{exackkk00}--\eqref{exackkk} guarantees the differential inequality
\vspace{+0mm}
\begin{equation}
\dfrac{d}{dt}V(x,\hat{\theta}) \leq -\mu x_{1}^{2}-\mu Q \left( x_{2}+\hat{\theta}_{1}x_{1}^{2}+\hat{\theta}_{2}x_{1}^{3}+\mu x_{1} \right)^{2}  \label{exdiein}
\end{equation}
for the Lyapunov function 
\vspace{+0mm}
\begin{align}
V(x,\hat{\theta}) =& \dfrac{1}{2}x_{1}^{2}+\dfrac{Q}{2} \left( x_{2}+\hat{\theta}_{1}x_{1}^{2}+\hat{\theta}_{2}x_{1}^{3}+\mu x_{1} \right)^{2}+\dfrac{1}{2\gamma_{1}} \left( \hat{\theta}_{1} - \theta_{1} \right)^{2}+\dfrac{1}{2\gamma_{2}} \left( \hat{\theta}_{2} - \theta_{2} \right)^{2} \label{exVdef} 
\end{align}
As shown in \cite{b1} the feedback law \eqref{exackkk00}--\eqref{exackkk} guarantees the following properties:

\begin{enumerate}
\item[(i)]the set $\left\lbrace (0,\hat{\theta}):\hat{\theta}\in \mathbb{R}^{p} \right\rbrace$ is the set of equilibrium points for the closed-loop system \eqref{eq:excs0}, \eqref{eq:excs} with \eqref{exackkk00}--\eqref{exackkk},
\item[(ii)]the equilibrium point $(x,\hat{\theta})=(0,\theta)$ is Lyapunov stable for the closed-loop system \eqref{eq:excs0}, \eqref{eq:excs} with \eqref{exackkk00}--\eqref{exackkk},
\item[(iii)]for every initial condition the solution $(x(t),\hat{\theta}(t))$ is bounded and system \eqref{eq:excs0}, \eqref{eq:excs} with \eqref{exackkk00}--\eqref{exackkk} is Lagrange stable,
\item[(iv)]the closed-loop system \eqref{eq:excs0}, \eqref{eq:excs} with \eqref{exackkk00}--\eqref{exackkk} and output $Y=x$ is GAOS and every solution satisfies $\lim_{t \rightarrow \infty}(x(t))=0$.
\end{enumerate}

If we follow the procedure described in the proof of Lemma \ref{lem4} and the proof of Theorem \ref{thm:2} then we obtain a different family of adaptive feedback controllers. More specifically, we get the adaptive feedback controller
\begin{align}
u =& -Q^{-1}x_{1}-(x_{1}+\hat{\theta}_{1})x_{1}w_{1}-(x_{1}^{2}+\hat{\theta}_{2})x_{1}w_{2}-Gz+\varphi(Mx_{1}-z) \label{eq:exmac00} 
\end{align}
\begin{align}
\dfrac{d\hat{\theta}_{1}}{dt} =& w_{1} =\gamma_{1}x_{1}^{2}(Qz\varphi+x_{1}) \label{eq:exmac0} \\
\dfrac{d\hat{\theta}_{2}}{dt} =& w_{2} =\gamma_{2}x_{1}^{3}(Qz\varphi+x_{1}) \label{eq:exmac}
\end{align}
where $Q, \gamma_{1} , \gamma_{2}, \mu, \epsilon >0$, $r \geq 0$ are constants and

\begin{align}
M :=& 2 \mu + \frac{\hat{\theta}_{1}^{2}+\hat{\theta}_{2}^{2}+r}{2} +x_{1}^{2}+\left( 1+\frac{1}{2\epsilon} \right) x_{1}^{4} +\frac{1}{2\epsilon}x_{1}^{6} \label{eq:exsmbdefmac000}  \\
z :=& x_{2} + \hat{\theta}_{1} x_{1}^{2} + \hat{\theta}_{2} x_{1}^{3} + Mx_{1} \label{eq:exsmbdefmac00} \\
\varphi :=& 2\hat{\theta}_{1}x_{1}+3\hat{\theta}_{2}x_{1}^{2}+ M +x_{1}^{2} \left( 2+2\frac{2\epsilon +1}{\epsilon} x_{1}^{2} +\frac{3}{\epsilon}x_{1}^{4} \right) \label{eq:exsmbdefmac0} \\
G :=& \mu +Qx_{1}^{2} \left( 1+x_{1}^{2} \right) \varphi^{2} \left( \frac{x_{1}^{2}}{2\epsilon}+\frac{\hat{\theta}_{1}^{2}+\hat{\theta}_{2}^{2}+r}{\mu} \right) \label{eq:exsmbdefmac}  
\end{align}

The adaptive controller \eqref{eq:exmac00}--\eqref{eq:exmac}, guarantees the following two differential inequalities:
\vspace{-0mm}
\begin{align}
\dfrac{d}{dt}W(x,\hat{\theta}) &\leq - \mu U(x,\hat{\theta}) \\
\dfrac{d}{dt}U(x,\hat{\theta}) &\leq - \mu U(x,\hat{\theta})+ \epsilon ( \vert \theta \vert^{2}-r)^{+}
\end{align}
\newpage

\begin{figure}[!h]
\centerline{\includegraphics[width=8.7cm]{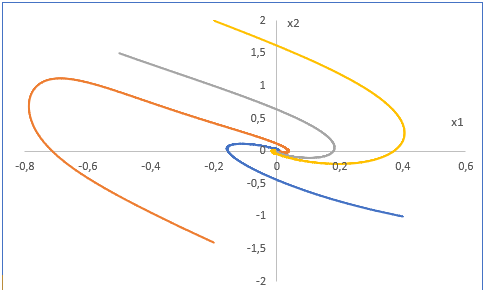}}
\caption{Projection of the $x_{1}-x_{2}$ plane of certain solutions of the closed-loop system \eqref{eq:excs0}, \eqref{eq:excs} with \eqref{exackkk00}--\eqref{exackkk}, initial condition $\left( \hat{\theta}_{1}(0),\hat{\theta}_{2}(0)\right)=(0,0)$, $\gamma_{1}=\gamma_{2}=\mu=Q=1$.}
\label{fig1}
\end{figure}

\vspace{+0mm}

\begin{figure}[!h]
\centerline{\includegraphics[width=8.7cm]{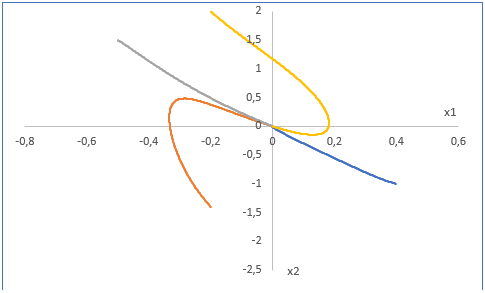}}
\caption{Projection of the $x_{1}-x_{2}$ plane of certain solutions of the closed-loop system \eqref{eq:excs0}, \eqref{eq:excs} with \eqref{eq:exmac00}--\eqref{eq:exmac}, initial condition $\left( \hat{\theta}_{1}(0),\hat{\theta}_{2}(0)\right)=(0,0)$, $\gamma_{1}=\gamma_{2}=\mu=Q=\epsilon=1$, $r=2$.}
\label{fig2}
\end{figure}

\vspace{+0mm}

\begin{figure}[!h]
\centerline{\includegraphics[width=8.7cm]{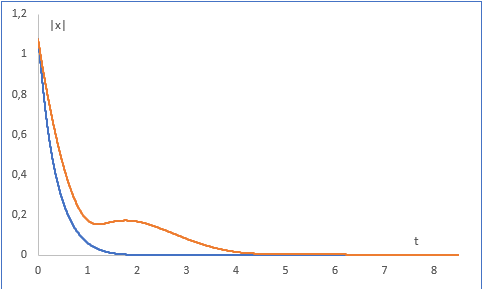}}
\caption{The evolution of $\vert x(t) \vert$ for the closed-loop system \eqref{eq:excs0}, \eqref{eq:excs} with \eqref{exackkk00}--\eqref{exackkk} (red line) and the closed-loop system \eqref{eq:excs} with \eqref{eq:exmac00}--\eqref{eq:exmac} (blue line), initial condition $\left( x_{1}(0),x_{2}(0),\hat{\theta}_{1}(0),\hat{\theta}_{2}(0)\right)=(.4,-1,0,0)$, $\gamma_{1}=\gamma_{2}=\mu=Q=\epsilon=1$, $r=2$}
\label{fig3}
\end{figure}

\begin{figure}[!h]
\centerline{\includegraphics[width=8.7cm]{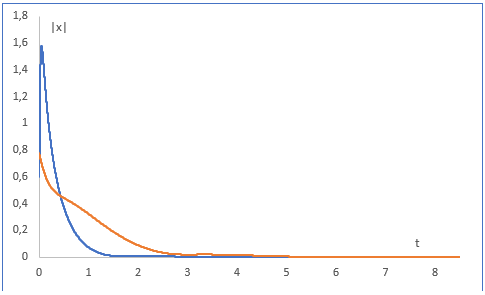}}
\caption{The evolution of $\vert x(t) \vert$ for the closed-loop system \eqref{eq:excs0}, \eqref{eq:excs} with \eqref{exackkk00}--\eqref{exackkk} (red line) and the closed-loop system \eqref{eq:excs0}, \eqref{eq:excs} with \eqref{eq:exmac00}--\eqref{eq:exmac} (blue line), initial condition $\left( x_{1}(0),x_{2}(0),\hat{\theta}_{1}(0),\hat{\theta}_{2}(0)\right)=(.6,.5,0,0)$, $\gamma_{1}=\gamma_{2}=\mu=Q=\epsilon=1$, $r=2$}
\label{fig4}
\end{figure}

\vspace{-0mm}

\begin{figure}[!h]
\centerline{\includegraphics[width=8.7cm]{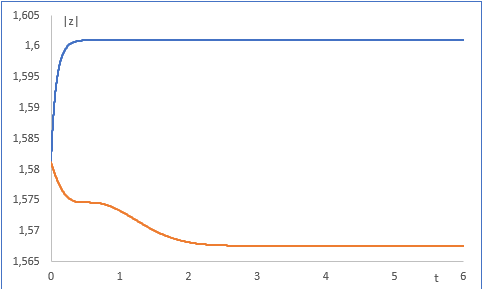}}
\caption{The evolution of $\vert z(t) \vert=\vert \hat{\theta}(t)-\theta \vert$ for the closed-loop system \eqref{eq:excs0}, \eqref{eq:excs} with \eqref{exackkk00}--\eqref{exackkk} (red line) and the closed-loop system \eqref{eq:excs0}, \eqref{eq:excs} with \eqref{eq:exmac00}--\eqref{eq:exmac} (blue line), initial condition $\left( x_{1}(0),x_{2}(0),\hat{\theta}_{1}(0),\hat{\theta}_{2}(0)\right)=(.4,-1,0,0)$, $\gamma_{1}=\gamma_{2}=\mu=Q=\epsilon=1$, $r=2$}
\label{fig5}
\end{figure}

\vspace{-0mm}

\begin{figure}[!h]
\centerline{\includegraphics[width=8.7cm]{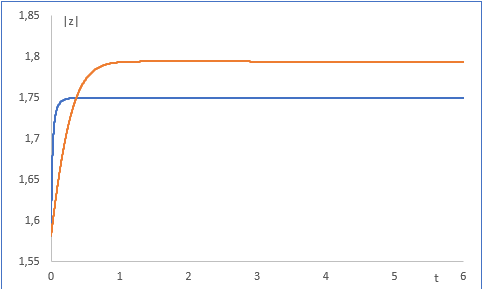}}
\caption{The evolution of $\vert z(t) \vert=\vert \hat{\theta}(t)-\theta \vert$ for the closed-loop system \eqref{eq:excs0}, \eqref{eq:excs} with \eqref{exackkk00}--\eqref{exackkk} (red line) and the closed-loop system \eqref{eq:excs0}, \eqref{eq:excs} with \eqref{eq:exmac00}--\eqref{eq:exmac} (blue line), initial condition $\left( x_{1}(0),x_{2}(0),\hat{\theta}_{1}(0),\hat{\theta}_{2}(0)\right)=(.6,.5,0,0)$, $\gamma_{1}=\gamma_{2}=\mu=Q=\epsilon=1$, $r=2$}
\label{fig6}
\end{figure}

\newpage
where 
\vspace{-0mm}
\begin{align}
U(x,\hat{\theta}) &:= \dfrac{1}{2}x_{1}^{2} + \dfrac{Q}{2}z^{2} \\
W(x,\hat{\theta}) &:= U(x,\hat{\theta}) + \dfrac{1}{2\gamma_{1}} \left( \hat{\theta}_{1} - \theta_{1} \right)^{2} + \dfrac{1}{2\gamma_{2}} \left( \hat{\theta}_{2} - \theta_{2} \right)^{2}
\end{align}
The difference between \eqref{exackkk00}--\eqref{exackkk} and \eqref{eq:exmac00}--\eqref{eq:exmac} is clearly visible. Theorem \ref{thm:2} guarantees the same properties as the feedback law \eqref{exackkk00}--\eqref{exackkk} and the additional following properties:

\begin{enumerate}
\item[(v)] the closed-loop system \eqref{eq:excs0}, \eqref{eq:excs} with \eqref{eq:exmac00}--\eqref{eq:exmac} and output $Y=x$ is UGAOS when $\vert \theta \vert \leq \sqrt{r} $,
\item[(vi)] the closed-loop system \eqref{eq:excs0}, \eqref{eq:excs} with \eqref{eq:exmac00}--\eqref{eq:exmac} is UGAOS for every output $Y=h(x,\hat{\theta},\theta)$ for which there exists a function $b \in K_{\infty}$ with \\ $\vert h(x,\hat{\theta},\theta) \vert \leq b\left( \left( x_{1}^{2}+z^{2}-\mu^{-1} \epsilon \left( \vert \theta \vert^{2}-r \right)^{+} \right)^{+} \right)$
\end{enumerate}
Although the difference between the behaviours of the two adaptive controllers is a qualitative difference and not a quantitative one, the projections of the solutions on the $x_{1}-x_{2}$ plane illustrate the qualitative change that takes place when the controller \eqref{eq:exmac00}--\eqref{eq:exmac} is used. This is shown in Figure \ref{fig1} and Figure \ref{fig2}.

The solutions of the closed-loop system \eqref{eq:excs0}, \eqref{eq:excs} with \eqref{eq:exmac00}--\eqref{eq:exmac} converge with a faster rate than the solutions of the closed-loop system \eqref{eq:excs0}, \eqref{eq:excs} with \eqref{exackkk00}--\eqref{exackkk} for the same parameter values. However, the solutions of the closed-loop system \eqref{eq:excs0}, \eqref{eq:excs} with \eqref{eq:exmac00}--\eqref{eq:exmac} may present a higher overshoot. This is shown in Figure \ref{fig3} and Figure \ref{fig4}. 

Moreover, for both closed-loop systems the parameter estimates $\hat{\theta}(t)$ do not converge to the real values of the parameter $\theta=\left( -\frac{3}{2},-\frac{1}{2} \right)$. This is shown in Figure \ref{fig5} and Figure \ref{fig6}. 

Figure \ref{fig5} also shows that the controller \eqref{eq:exmac00}--\eqref{eq:exmac} exhibits less effort for the estimation of the unknown parameters than the controller \eqref{exackkk00}--\eqref{exackkk}. Indeed, the controller \eqref{eq:exmac00}--\eqref{eq:exmac} acts primarily as a robust nonlinear controller and thus the parameter estimation becomes less important.

\section{\text{Proofs}}

\subsection{\text{Proof} of Theorem \ref{thm:1}}
We use definition \eqref{eq:mac}, equation \eqref{erl}, the Young inequalities:
\vspace{-0mm}
\begin{align}
\nabla P(x) & g(x)(\varphi(x))'\left( \theta- \bar{P}(\theta)\right) \leq  \frac{1}{2\delta} \vert \theta- \bar{P}(\theta)\vert^{2}+\frac{\delta}{2}( \nabla P(x)g(x))^{2}\vert \varphi(x)\vert^{2}
\end{align}

\begin{align}
-\nabla P(x) & g(x)(\varphi(x))' \hat{\theta} \leq  \mu (x)( \nabla P(x) g(x))^{2} \vert \hat{\theta} \vert^{2} +\frac{1}{4\mu(x)} \vert \varphi(x)\vert^{2} 
\end{align}

\begin{align} 
\nabla P(x) & g(x)(\varphi(x))'\bar{P}(\theta) \leq  \mu (x)( \nabla P(x) g(x))^{2} \vert \bar{P}(\theta) \vert^{2} +\frac{1}{4\mu(x)} \vert \varphi(x)\vert^{2}
\end{align}
where $\bar{P}:\mathbb{R}^{p} \rightarrow \mathbb{R}^{p}$ the projection of $\theta$ on the ball of radius $\sqrt{r}$ centered at $0 \in \mathbb{R}^{p}$

\begin{equation}
\bar{P}(\theta):= \begin{cases} 
\theta, &\text{ if } \vert \theta \vert \leq \sqrt{r} \\
\frac{\theta}{\vert  \theta \vert} \sqrt{r} &\text{ if } \vert \theta \vert > \sqrt{r}
\end{cases} \label{eq:projdef}
\end{equation}
and inequality \eqref{hh} in order to obtain the estimate:
\vspace{-0mm}
\begin{align}
\frac{d}{dt}P(x) \leq& -Q(x) +\frac{1}{2\mu (x)} \vert \varphi(x)\vert^{2}+ \frac{1}{2\delta} \vert \theta- \bar{P}(\theta)\vert^{2}+\mu (x) (\vert \bar{P}(\theta) \vert^{2}-r)(\nabla P(x) g(x))^{2} \label{eq:th1pr1in}
\end{align}
for all $(x,\hat{\theta}) \in \mathbb{R}^{n} \times \mathbb{R}^{p}$.

Using inequalities \eqref{aa}, \eqref{eq:t1h}, \eqref{eq:th1pr1in} and the facts that $ \vert \bar{P}(\theta) \vert^{2} \leq r $, $\vert \theta -\bar{P}(\theta) \vert^{2} \leq \left( \vert \theta \vert^{2} -r  \right)^{+} $ (consequences of definition \eqref{eq:projdef}), we obtain:

\begin{equation}
\frac{d}{dt}P(x) \leq -\frac{1}{2}\rho(P(x)) + \frac{1}{2\delta} \left( \vert \theta \vert^{2} -r  \right)^{+} \label{eq:th1pr2in}
\end{equation}
for all $(x,\hat{\theta}) \in \mathbb{R}^{n} \times \mathbb{R}^{p}$.
We define the function 

\begin{equation}
W(x)=\frac{1}{2}\left( \left( P(x)-\alpha \right)^{+} \right)^{2}, \text{ for } x \in \mathbb{R}^{n} \label{eq:th1prdefW}
\end{equation}
and we notice that inequality \eqref{eq:th1pr2in} in conjunction with definition \eqref{eq:defa} guarantees the following implication:

\begin{equation}
W(x)>0 \Rightarrow \frac{d}{dt}P(x) \leq -\frac{\lambda}{2} \rho(P(x)) \label{eq:th1primpl}
\end{equation}
Since $\frac{d}{dt} W(x)=\left( P(x)-\alpha \right)^{+} \frac{d}{dt} P(x)$, using \eqref{eq:th1primpl} and distinguishing the cases $W(x)=0$ and $W(x)>0$, we get for all $(x,\hat{\theta}) \in \mathbb{R}^{n} \times \mathbb{R}^{p}$:

\begin{align}
\frac{d}{dt} W(x) &\leq -\frac{\lambda}{2} \left( P(x)-\alpha \right)^{+} \rho(P(x)) \nonumber \\ 
& \leq -\frac{\lambda}{2} \left( P(x)-\alpha \right)^{+} \rho(\left( P(x)-\alpha \right)^{+}) \nonumber \\
& = -\tilde{\rho}(W(x))
\end{align}
where $\tilde{\rho}(s):=\lambda\sqrt{\frac{s}{2}}\rho(\sqrt{2s})$ is a function of class $K_{\infty}$. The existence of $\sigma \in KL$ for which \eqref{eq:gkl1} holds is a direct consequence of the above differential inequality and Lemma 4.4 in \cite{b4}.

\subsection{Proof of Lemma 2}

We have $h(x,\theta)=\sum_{i=1}^{n}\int_{0}^{1}\dfrac{\partial h}{\partial x_{i}}(\lambda x,\theta)x_{i}d\lambda$ for all $(x,\theta) \in \mathbb{R}^{n} \times \mathbb{R}^{p}$. Using the Cauchy-Schwarz inequality, we get 
\begin{equation}
\vert h(x,\theta) \vert \leq \left( \sum_{i=1}^{n} \left( \int_{0}^{1}\dfrac{\partial h}{\partial x_{i}}(\lambda x,\theta)d\lambda \right)^{2}  \right)^{\frac{1}{2}} \vert x \vert
\end{equation}
for all $(x,\theta) \in \mathbb{R}^{n} \times \mathbb{R}^{p}$. Again, using the 
Cauchy-Schwarz inequality, we get
\begin{equation}
\vert h(x,\theta) \vert \leq \left( \sum_{i=1}^{n}  \int_{0}^{1} \left( \dfrac{\partial h}{\partial x_{i}}(\lambda x,\theta) \right)^{2} d\lambda   \right)^{\frac{1}{2}} \vert x \vert
\end{equation}
for all $(x,\theta) \in \mathbb{R}^{n} \times \mathbb{R}^{p}$. We conclude that the inequality $\vert h(x,\theta) \vert \leq \rho(x,\theta) \vert x \vert$ holds for all $(x,\theta) \in \mathbb{R}^{n} \times \mathbb{R}^{p}$ with 
\begin{equation}
\rho(x,\theta):=\left( 1+ \sum_{i=1}^{n}  \int_{0}^{1} \left( \dfrac{\partial h}{\partial x_{i}}(\lambda x,\theta) \right)^{2} d\lambda   \right)^{\frac{1}{2}}
\end{equation} 
\begin{flushright}
$\blacksquare$
\end{flushright}
\vspace{-0mm}
\subsection{Proof of Lemma 3}

Notice that $h(x,\theta)=h(T^{-1}(T(x,\theta),\theta),\theta)$ for all $(x,\theta) \in \mathbb{R}^{n} \times \mathbb{R}^{p}$. Consider the map $\tilde{h}=h(T^{-1}(z,\theta),\theta)$ which is a $C^{\infty}$ map that satisfies $\tilde{h}(0,\theta)=0$ for all $\theta \in \mathbb{R}^{p}$. By virtue of Lemma \ref{lem1}, there exists a $C^{\infty}$ map $\tilde{\rho}:\mathbb{R}^{n} \times \mathbb{R}^{p} \rightarrow \left[1, +\infty \right) $ such that  $\vert \tilde{h}(z,\theta) \vert \leq \tilde{\rho}(z,\theta) \vert z \vert$ for all $(z,\theta) \in \mathbb{R}^{n} \times \mathbb{R}^{p}$. Setting $z=T(x,\theta)$ and noticing that $h(x,\theta)=\tilde{h}(T(x,\theta),\theta)$ we get $\vert h(x,\theta) \vert \leq \tilde{\rho}(T(x,\theta),\theta) \vert T(x,\theta) \vert$ for all $(x,\theta) \in \mathbb{R}^{n} \times \mathbb{R}^{p}$. We conclude that the inequality $\vert h(x,\theta) \vert \leq \rho(x,\theta) \vert T(x,\theta) \vert$ holds for all $(x,\theta) \in \mathbb{R}^{n} \times \mathbb{R}^{p}$ with $\rho(x,\theta):= \tilde{\rho}(T(x,\theta),\theta).$
\begin{flushright}
$\blacksquare$
\end{flushright}
\vspace{-0mm}
\subsection{Proof of Lemma 4}

Let $\bar{\Theta}\subset \mathbb{R}^{p}$ be the closure of $\Theta\subset \mathbb{R}^{p}$ (a compact set). Define $b(s):= max\left\lbrace \vert T(x,\theta) \vert : x \in \mathbb{R}^{n}, \theta \in \bar{\Theta}, \vert x \vert \leq s \right\rbrace$ for $s \geq 0$. By definition $b :\mathbb{R}_{+} \rightarrow \mathbb{R}_{+}$ is a well defined, non-decreasing function with $b(0)=0$. Continuity of $T:\mathbb{R}^{n} \times \mathbb{R}^{p} \rightarrow \mathbb{R}^{p}$, the fact that $T(0,\theta)=0$ for all $\theta \in \mathbb{R}^{p}$ and compactness of $\bar{\Theta}\subset \mathbb{R}^{p}$ implies that for every $\epsilon>0$ there exists $\delta>0$ with $\vert T(x,\theta) \vert < \epsilon$ for all $(x,\theta) \in \mathbb{R}^{n} \times \Theta$ with $\vert x \vert < \delta$. Therefore, $\lim_{s\rightarrow 0^{+}}(b(s))=b(0)=0$. Lemma 2.4 in \cite{b3} implies the existence of a $K_{\infty}$ function $\beta :\mathbb{R}_{+} \rightarrow \mathbb{R}_{+}$ such that $b(s) \leq \beta(s)$ for all $s \geq 0$. Therefore, it holds that $\vert T(x,\theta) \vert \leq \beta(x)$ for all $(x,\theta) \in \mathbb{R}^{n} \times \Theta$. Similarly, by defining $c(s):=max \left\lbrace \vert T^{-1}(x,\theta) \vert : x \in \mathbb{R}^{n}, \theta \in \bar{\Theta}, \vert x \vert \leq s \right\rbrace$ for $s \geq 0$, we conclude the existence of a $K_{\infty}$ function $C:\mathbb{R}_{+} \rightarrow \mathbb{R}_{+}$ such that $c(s)\leq C(s)$ for all $s\geq 0$. Therefore, it holds that $\vert T^{-1}(x,\theta) \vert \leq C(\vert x \vert)$ for all $(x,\theta) \in \mathbb{R}^{n} \times \Theta$. The fact that $T^{-1}(T(x,\theta),\theta)=x$ for all $(x,\theta) \in \mathbb{R}^{n} \times \mathbb{R}^{p}$ implies that $\vert x \vert \leq C(\vert T(x,\theta) \vert)$ for all $(x,\theta) \in \mathbb{R}^{n} \times \Theta$. Consequently, it holds that $\vert T(x,\theta) \vert \geq \alpha (\vert x \vert)$ for all $(x,\theta) \in \mathbb{R}^{n} \times \Theta$ with $\alpha =C^{-1}$.
\begin{flushright}
$\blacksquare$
\end{flushright}
\vspace{-0mm}
\subsection{Proof of Lemma 5}

Define the map $\tilde{T}:\mathbb{R}^{n+1} \times \mathbb{R}^{p} \rightarrow \mathbb{R}^{n+1}$ for all $(x,y,\hat{\theta}) \in \mathbb{R}^{n+1} \times \mathbb{R}^{p}$ by means of the following equation
\vspace{+1mm}
\begin{equation}
\tilde{T}(x,y,\hat{\theta}):= \begin{pmatrix}
T(x,\hat{\theta})\\
y-k(x,\hat{\theta})
\end{pmatrix} \label{eq:pl4dd}
\end{equation}
Since $T(0,\hat{\theta})=0$ and $k(0,\hat{\theta})=0$ for all $\hat{\theta} \in \mathbb{R}^{p}$, it follows that \eqref{eq:pl4dd} that $T(0,0,\hat{\theta})=0$ for all $\hat{\theta} \in \mathbb{R}^{p}$. Since $T:\mathbb{R}^{n} \times \mathbb{R}^{p} \rightarrow \mathbb{R}^{n}$ and $k:\mathbb{R}^{n} \times \mathbb{R}^{p} \rightarrow \mathbb{R}^{n}$ are $C^{\infty}$ maps, it follows from \eqref{eq:pl4dd} that $\tilde{T}:\mathbb{R}^{n+1} \times \mathbb{R}^{p} \rightarrow \mathbb{R}^{n+1}$ is a $C^{\infty}$ map. The map $\tilde{T}^{-1}:\mathbb{R}^{n+1} \times \mathbb{R}^{p} \rightarrow \mathbb{R}^{n+1}$ defined for all $(z,\zeta,\hat{\theta}) \in \mathbb{R}^{n+1} \times \mathbb{R}^{p}$ by the equation
\vspace{+0mm} 
\begin{equation}
\tilde{T}^{-1}(z,\zeta,\hat{\theta}):= \begin{pmatrix}
T^{-1}(z,\hat{\theta})\\
\zeta + k(T^{-1}(z,\hat{\theta}),\hat{\theta}) \label{eq:pl4did}
\end{pmatrix}
\end{equation}
is also a $C^{\infty}$ map that satisfies $\tilde{T}(\tilde{T}^{-1}(x,y,\hat{\theta}),\hat{\theta})=\tilde{T}^{-1}(\tilde{T}(x,y,\hat{\theta}),\hat{\theta})=(x,y)'$ for all  $(x,y,\hat{\theta}) \in \mathbb{R}^{n+1} \times \mathbb{R}^{p}$. Hence, the map $\tilde{T}:\mathbb{R}^{n+1} \times \mathbb{R}^{p} \rightarrow \mathbb{R}^{n+1}$ defined by \eqref{eq:pl4dd} is a smooth parameterized family of diffeomorphisms of $\mathbb{R}^{n+1}$.
Next, define the update law $\tilde{w}(x,y,\hat{\theta}) = \left( \tilde{w}_{1}(x,y,\hat{\theta}),\ldots,\tilde{w}_{p}(x,y,\hat{\theta}) \right)' \in \mathbb{R}^{p}$ for all $(x,y,\hat{\theta}) \in \mathbb{R}^{n+1} \times \mathbb{R}^{p}$ by means of the following formula for $j=1,\ldots,p$:
\begin{align}
\tilde{w}_{j}(x,y,\hat{\theta}):=& w_{j}(x,\hat{\theta})+\gamma_{j}(y-k(x,\hat{\theta}))  \left( \varphi_{n+1,j}(x,y)-\sum_{i=1}^{n}\dfrac{\partial k}{\partial x_{i}}(x,\hat{\theta}) \varphi_{i,j}(x) \right) \label{eq:pl4uld}
\end{align}
\vspace{+0mm}
It should be noticed that inequalities \eqref{eq:inV}, \eqref{eq:inrt} and definition \eqref{eq:defV} imply the following inequalities for all $(x,\hat{\theta},\theta) \in \mathbb{R}^{n} \times \mathbb{R}^{p} \times \mathbb{R}^{p}$:
\vspace{-0mm} 
\begin{align}
& \sum_{i=1}^{n}(T(x,\hat{\theta}))'\dfrac{\partial T}{\partial x_{i}}(x,\hat{\theta})\left( F_{i}(x)+b_{i}(x)k(x,\hat{\theta})+\sum_{j=1}^{p}\varphi_{i,j}\theta_{j} \right) \nonumber \\
& +\sum_{j=1}^{p}(T(x,\hat{\theta}))'\dfrac{\partial T}{\partial \hat{\theta}_{j}}(x,\hat{\theta})w_{j}(x,\hat{\theta})  \nonumber \\  \label{eq:48}
& +\sum_{j=1}^{p}\frac{1}{\gamma_{j}}(\hat{\theta}_{j}-\theta_{j})w_{j}(x,\hat{\theta}) \leq -\alpha \vert T(x,\hat{\theta}) \vert^{2}  
\end{align} 
\begin{align}
& \sum_{i=1}^{n}(T(x,\hat{\theta}))'\dfrac{\partial T}{\partial x_{i}}(x,\hat{\theta})\left( F_{i}(x)+b_{i}(x)k(x,\hat{\theta})+\sum_{j=1}^{p}\varphi_{i,j}\theta_{j} \right) \nonumber \\
& +\sum_{j=1}^{p}(T(x,\hat{\theta}))'\dfrac{\partial T}{\partial \hat{\theta}_{j}}(x,\hat{\theta})w_{j}(x,\hat{\theta})  \nonumber \\
& \leq -\omega \vert T(x,\hat{\theta}) \vert^{2}+ \epsilon \left( \vert \theta \vert^{2}-r \right)^{+}   \label{eq:47}
\end{align}
\vspace{-0mm}
We next define the feedback law $\tilde{k}(x,y,\hat{\theta})$ for all $(x,y,\hat{\theta}) \in \mathbb{R}^{n+1} \times \mathbb{R}^{p}$ by means of the following formula:
\vspace{-0mm}
\begin{align} 
\tilde{k}(x,y,\hat{\theta}):=& \dfrac{\bar{k}(x,y,\hat{\theta})}{g(x,y)}  \\
\bar{k}(x,y,\hat{\theta}):=& -f(x,y)- \sum_{j=1}^{p}\varphi_{i,j}(x,y)\hat{\theta}_{j}+ \sum_{i=1}^{n}\dfrac{\partial k}{\partial x_{i}}(x,\hat{\theta}) \left( F_{i}(x)+b_{i}(x)y + \sum_{j=1}^{p}\varphi_{i,j}(x)\hat{\theta}_{j} \right) \nonumber \\
& +\sum_{j=1}^{p}\dfrac{\partial k}{\partial \hat{\theta}_{j}}(x,\hat{\theta})w_{j}(x,\hat{\theta}) +\sum_{j=1}^{p}\dfrac{\partial k}{\partial \hat{\theta}_{j}}(x,\hat{\theta})\gamma_{j}(y-k(x,\hat{\theta}))\left( \varphi_{n+1,j}(x,y)-\sum_{i=1}^{n}\dfrac{\partial k}{\partial x_{i}}(x,\hat{\theta})\varphi_{i,j}(x) \right) \nonumber \\
& -\sum_{j=1}^{p}(T(x,\hat{\theta}))'\dfrac{\partial T}{\partial \hat{\theta}_{j}}(x,\hat{\theta})\gamma_{j}  \left( \varphi_{n+1,j}(x,y)-\sum_{i=1}^{n}\dfrac{\partial k}{\partial x_{i}}(x,\hat{\theta})\varphi_{i,j}(x) \right) \nonumber \\  
& -\sum_{i=1}^{n}b_{i}(T(x,\hat{\theta}))'\dfrac{\partial T}{\partial x_{i}}(x,\hat{\theta})- M(x,y,\hat{\theta})(y-k(x,\hat{\theta}))  \label{eq:49}  
\end{align}
\vspace{-0mm}
where $M(x,y,\hat{\theta})$ is a positive $C^{\infty}$ function that is to be determined next.  
Notice that by virtue of definition \eqref{eq:pl4dd}, the closed-loop system \eqref{eq:cslemma}, \eqref{eq:2eqcslem} with the adaptive feedback law \eqref{eq:afl2lem0}, \eqref{eq:afl2lem} satisfies the following equation:
\vspace{-0mm}
\begin{align}
 \frac{d}{dt} \left( \frac{1}{2} \vert \tilde{T}(x,y,\hat{\theta}) \vert^{2} \right)=&  \sum_{i=1}^{n}(T(x,\hat{\theta}))'\dfrac{\partial T}{\partial x_{i}}(x,\hat{\theta})\left( F_{i}(x)+b_{i}(x)y+\sum_{j=1}^{p}\varphi_{i,j}\theta_{j} \right) \nonumber \\
&+\sum_{j=1}^{p}(T(x,\hat{\theta}))'\dfrac{\partial T}{\partial \hat{\theta}_{j}}(x,\hat{\theta})\tilde{w}(x,y,\hat{\theta})- (y-k(x,\hat{\theta}))\sum_{j=1}^{p}\dfrac{\partial k}{\partial \hat{\theta}_{j}}(x,\hat{\theta})\tilde{w}_{j}(x,y,\hat{\theta}) \nonumber \\
&+(y-k(x,\hat{\theta}))  \left( f(x,y)+g(x,y)\tilde{k}(x,y,\hat{\theta})+\sum_{j=1}^{p}\varphi_{n+1,j}(x,y)\theta_{j} \right) \nonumber \\
& -(y-k(x,\hat{\theta})) \sum_{i=1}^{n}\dfrac{\partial k}{\partial x_{i}} \left( F_{i}(x)+b_{i}(x)y+\sum_{j=1}^{p}\varphi_{i,j}(x)\theta_{j} \right) \label{eq:pl4eq1}
\end{align}

and the equation:

\begin{align}
 \frac{d}{dt} \left( \frac{1}{2} \vert \tilde{T}(x,y,\hat{\theta}) \vert^{2} +\sum_{j=1}^{p}\frac{1}{2\gamma_{j}}\left( \hat{\theta}_{j}-\theta_{j} \right)^{2} \right) =& 
 \sum_{i=1}^{n}(T(x,\hat{\theta}))'\dfrac{\partial T}{\partial x_{i}}(x,\hat{\theta})\left( F_{i}(x)+b_{i}(x)y+\sum_{j=1}^{p}\varphi_{i,j}\theta_{j} \right) \nonumber \\
& +\sum_{j=1}^{p}(T(x,\hat{\theta}))'\dfrac{\partial T}{\partial \hat{\theta}_{j}}(x,\hat{\theta})\tilde{w}(x,y,\hat{\theta}) +\sum_{j=1}^{p}\frac{1}{\gamma_{j}}\left( \hat{\theta}_{j}-\theta_{j} \right)\tilde{w}(x,y,\hat{\theta}) \nonumber \\
& +(y-k(x,\hat{\theta}))\left( f(x,y)+g(x,y)\tilde{k}(x,y,\hat{\theta})+\right. 
\nonumber \\
&\left.  +\sum_{j=1}^{p}\varphi_{n+1,j}(x,y)\theta_{j} -\sum_{j=1}^{p}\frac{\partial k}{\partial \hat{\theta}_{j}}(x,\hat{\theta})\tilde{w}(x,y,\hat{\theta}) \right) \nonumber \\
& -(y-k(x,\hat{\theta}))  \sum_{i=1}^{n}\dfrac{\partial k}{\partial x_{i}} \left( F_{i}(x)+b_{i}(x)y+\sum_{j=1}^{p}\varphi_{i,j}(x)\theta_{j} \right) \label{eq:pl4eq2}
\end{align}
Combining \eqref{eq:pl4uld}, \eqref{eq:47}, \eqref{eq:48}, \eqref{eq:49}, \eqref{eq:pl4eq1}, \eqref{eq:pl4eq2} and definition \eqref{eq:defV2} we get for all $(x,y,\hat{\theta}) \in \mathbb{R}^{n+1} \times \mathbb{R}^{p}$:

\begin{align}
\frac{d}{dt}\tilde{V}(x,y,\hat{\theta}) \leq & -\alpha \vert T(x,\hat{\theta}) \vert^{2}  - M(x,y,\hat{\theta}) \left( y-k(x,\hat{\theta}) \right)^{2} \label{eq:pl419}
\end{align}

\begin{align}
 \frac{d}{dt}\left( \frac{1}{2}\vert \tilde{T}(x,y,\hat{\theta}) \vert^{2} \right) \leq& -\omega \vert T(x,\hat{\theta}) \vert^{2} +\epsilon \left( \vert \theta \vert^{2}-r \right)^{+}  - M(x,y,\hat{\theta}) \left( y-k(x,\hat{\theta}) \right)^{2} \nonumber \\
& + \left( y-k(x,\hat{\theta}) \right)   \sum_{j=1}^{p}\left( \varphi_{n+1,j}(x,y)-\sum_{i=1}^{n}\frac{\partial k}{\partial x_{i}}(x,\hat{\theta})\varphi_{i,j}(x) \right) \left( \hat{\theta}_{j}-\theta_{j} \right) \label{eq:pl420}
\end{align}
By virtue of Lemma \ref{lem2}, there exist $C^{\infty}$ maps $\rho_{j}:\mathbb{R}^{n+1} \times \mathbb{R}^{p} \rightarrow \left[1,+\infty \right) $, $j=1,\ldots,p$ for which the following inequalities hold for $j=1,\ldots,p$:

\begin{equation}
\vert \varphi_{n+1,j}(x,y)-\sum_{i=1}^{n}\frac{\partial k}{\partial x_{i}}(x,\hat{\theta})\varphi_{i,j}(x) \vert \leq \rho_{j}(x,y,\hat{\theta}) \vert \tilde{T}(x,y,\hat{\theta}) \vert \label{eq:pl421}
\end{equation}
for all $(x,y,\hat{\theta}) \in \mathbb{R}^{n+1} \times \mathbb{R}^{p}$.

Definition \eqref{eq:pl4dd} in conjunction with inequalities \eqref{eq:pl421} gives the following inequalities for $j=1,\ldots,p$:

\begin{align}
\vert \varphi_{n+1,j}&(x,y) -\sum_{i=1}^{n}\frac{\partial k}{\partial x_{i}}(x,\hat{\theta})\varphi_{i,j}(x) \vert \leq \rho_{j}(x,y,\hat{\theta}) \vert T(x,\hat{\theta}) \vert + \rho_{j}(x,y,\hat{\theta}) \vert y-k(x,\hat{\theta}) \vert \label{eq:pl422}
\end{align}
for all $(x,y,\hat{\theta}) \in \mathbb{R}^{n+1} \times \mathbb{R}^{p}$.

Define the function $\bar{P}:\mathbb{R}^{p}\rightarrow\mathbb{R}^{p}$ as in \eqref{eq:projdef}. We get from \eqref{eq:pl420} for all $(x,y,\hat{\theta}) \in \mathbb{R}^{n+1} \times \mathbb{R}^{p}$:

\begin{align}
 \frac{d}{dt}\left( \frac{1}{2}\vert \tilde{T}(x,y,\hat{\theta}) \vert^{2} \right) \leq& -\omega \vert T(x,\hat{\theta}) \vert^{2} +\epsilon \left( \vert \theta \vert^{2}-r \right)^{+}  - M(x,y,\hat{\theta}) \left( y-k(x,\hat{\theta}) \right)^{2} \nonumber \\
& +\left( y-k(x,\hat{\theta}) \right)  \sum_{j=1}^{p}\left( \varphi_{n+1,j}(x,y)-\sum_{i=1}^{n}\frac{\partial k}{\partial x_{i}}(x,\hat{\theta})\varphi_{i,j}(x) \right) \bar{P}_{j}(\theta) \nonumber \\
& +\left( y-k(x,\hat{\theta}) \right)  \sum_{j=1}^{p}\left( \varphi_{n+1,j}(x,y)-\sum_{i=1}^{n}\frac{\partial k}{\partial x_{i}}(x,\hat{\theta})\varphi_{i,j}(x) \right) \left( \theta_{j}-\bar{P}_{j}(\theta) \right) \nonumber \\
& -\left( y-k(x,\hat{\theta}) \right)  \sum_{j=1}^{p}\left( \varphi_{n+1,j}(x,y)-\sum_{i=1}^{n}\frac{\partial k}{\partial x_{i}}(x,\hat{\theta})\varphi_{i,j}(x) \right) \hat{\theta}_{j} \label{eq:pl424}
\end{align}  
The inequalities 
\begin{align}
& \left( y-k(x,\hat{\theta}) \right) \sum_{j=1}^{p}\left( \varphi_{n+1,j}(x,y)-\sum_{i=1}^{n}\frac{\partial k}{\partial x_{i}}(x,\hat{\theta})\varphi_{i,j}(x) \right)   \nonumber \\
& \times \left( \theta_{j}-\bar{P}_{j}(\theta) \right) \leq \epsilon \left( \theta_{j}-\bar{P}_{j}(\theta) \right)^{2}+\frac{1}{4\epsilon}\left( y-k(x,\hat{\theta}) \right)^{2} \nonumber \\
& \times \left( \varphi_{n+1,j}(x,y)-\sum_{i=1}^{n}\frac{\partial k}{\partial x_{i}}(x,\hat{\theta})\varphi_{i,j}(x) \right)^{2} \label{eq:pl424,5}
\end{align}
which hold for all $j=1,\ldots,p$ and the inequalities \eqref{eq:pl422} combined with \eqref{eq:pl424} give for all $(x,y,\hat{\theta}) \in \mathbb{R}^{n+1} \times \mathbb{R}^{p}$:

\begin{align}
\frac{d}{dt}\left( \frac{1}{2}\vert \tilde{T}(x,y,\hat{\theta}) \vert^{2} \right) \leq& -\omega \vert T(x,\hat{\theta}) \vert^{2} +\epsilon \left( \vert \theta \vert^{2}-r \right)^{+} - M(x,y,\hat{\theta}) \left( y-k(x,\hat{\theta}) \right)^{2} \nonumber \\
& + \vert y-k(x,\hat{\theta}) \vert \vert T(x,\hat{\theta})\vert \sum_{j=1}^{p}\rho(x,y,\hat{\theta}) \vert \bar{P}_{j}(\theta) \vert+\epsilon \sum_{j=1}^{p}\left( \theta_{j}-\bar{P}_{j}(\theta) \right)^{2} \nonumber \\
& + \left( y-k(x,\hat{\theta}) \right)^{2} \sum_{j=1}^{p} \left( \vert \bar{P}_{j}(\theta) \vert+\vert \hat{\theta}_{j} \vert \right)\rho(x,y,\hat{\theta}) \nonumber \\
& +\frac{1}{4\epsilon}\left( y-k(x,\hat{\theta}) \right)^{2}  \left( \varphi_{n+1,j}(x,y)-\sum_{i=1}^{n}\frac{\partial k}{\partial x_{i}}(x,\hat{\theta})\varphi_{i,j}(x) \right)^{2} \nonumber \\
& + \vert y-k(x,\hat{\theta}) \vert \vert T(x,\hat{\theta}) \vert \sum_{j=1}^{p}\rho(x,y,\hat{\theta}) \vert \hat{\theta}_{j} \vert \label{eq:pl425} 
\end{align} 
Definition \eqref{eq:projdef} implies that $\vert \bar{P}_{j}(\theta) \vert \leq \sqrt{r}$ for $j=1,\ldots,p$ and $\sum_{j=1}^{p}\left( \theta_{j}-\bar{P}_{j}(\theta) \right)^{2} \leq \left( \vert \theta \vert^{2}-r \right)^{+}$. Using the previous inequalities in conjunction with \eqref{eq:pl425} and the inequalities  
\begin{align}
\vert \hat{\theta}_{j}\vert &\leq \vert \hat{\theta} \vert \\   
\vert \hat{\theta}_{j} \vert &\leq \frac{1}{2}+\frac{1}{2}\hat{\theta}_{j}^{2}
\end{align}
for $j=1,\ldots,p$, we obtain for all $(x,y,\hat{\theta}) \in \mathbb{R}^{n+1} \times \mathbb{R}^{p}$:

\begin{align}
\frac{d}{dt}\left( \frac{1}{2}\vert \tilde{T}(x,y,\hat{\theta}) \vert^{2} \right) \leq& -\omega \vert T(x,\hat{\theta}) \vert^{2} +2\epsilon \left( \vert \theta \vert^{2}-r \right)^{+}- M(x,y,\hat{\theta}) \left( y-k(x,\hat{\theta}) \right)^{2} \nonumber \\
& + \vert y-k(x,\hat{\theta}) \vert \vert T(x,\hat{\theta})\vert \sqrt{r} \sum_{j=1}^{p}\rho(x,y,\hat{\theta})  \nonumber \\
& +\left( y-k(x,\hat{\theta}) \right)^{2} \sum_{j=1}^{p} \left( \sqrt{r}+\frac{1}{2}+\frac{1}{2}\hat{\theta}_{j}^{2} \right) \rho(x,y,\hat{\theta}) \nonumber \\
& +\frac{1}{4\epsilon}\left( y-k(x,\hat{\theta}) \right)^{2}  \left( \varphi_{n+1,j}(x,y)-\sum_{i=1}^{n}\frac{\partial k}{\partial x_{i}}(x,\hat{\theta})\varphi_{i,j}(x) \right)^{2} \nonumber \\
& +\vert y-k(x,\hat{\theta}) \vert \vert T(x,\hat{\theta}) \vert \vert \hat{\theta} \vert \sum_{j=1}^{p}\rho(x,y,\hat{\theta})  \label{eq:pl426}
\end{align}
Using the inequalities        
\begin{align}
 \vert y-k(x,\hat{\theta}) \vert \vert T(x,\hat{\theta})\vert \sqrt{r} \sum_{j=1}^{p}\rho(x,y,\hat{\theta})  
 \leq& \frac{\omega}{4}\vert T(x,\hat{\theta})\vert^{2}+\frac{r}{\omega}\left( y-k(x,\hat{\theta}) \right)^{2} \left( \sum_{j=1}^{p}\rho(x,y,\hat{\theta}) \right)^{2} \\
 \vert y-k(x,\hat{\theta}) \vert \vert T(x,\hat{\theta}) \vert \vert \hat{\theta} \vert \sum_{j=1}^{p}\rho(x,y,\hat{\theta}) 
 \leq& \frac{\omega}{4}\vert T(x,\hat{\theta})\vert^{2}+\frac{\vert \hat{\theta} \vert}{\omega}\left( y-k(x,\hat{\theta}) \right)^{2}  \left( \sum_{j=1}^{p}\rho(x,y,\hat{\theta}) \right)^{2} 
\end{align}
we obtain from \eqref{eq:pl426} for all $(x,y,\hat{\theta}) \in \mathbb{R}^{n+1} \times \mathbb{R}^{p}$:
\begin{align}
\frac{d}{dt}\left( \frac{1}{2}\vert \tilde{T}(x,y,\hat{\theta}) \vert^{2} \right) \leq& -\omega \vert T(x,\hat{\theta}) \vert^{2} +2\epsilon \left( \vert \theta \vert^{2}-r \right)^{+}  - M(x,y,\hat{\theta}) \left( y-k(x,\hat{\theta}) \right)^{2} \nonumber \\
& +\frac{r+\vert \hat{\theta} \vert}{\omega}\left( y-k(x,\hat{\theta}) \right)^{2} \left( \sum_{j=1}^{p}\rho(x,y,\hat{\theta}) \right)^{2} \nonumber \\
& +\left( y-k(x,\hat{\theta}) \right)^{2} \sum_{j=1}^{p} \left( \sqrt{r}+\frac{1}{2}+\frac{1}{2}\hat{\theta}_{j}^{2} \right) \rho(x,y,\hat{\theta}) \nonumber \\
& +\frac{1}{4\epsilon}\left( y-k(x,\hat{\theta}) \right)^{2}  \left( \varphi_{n+1,j}(x,y)-\sum_{i=1}^{n}\frac{\partial k}{\partial x_{i}}(x,\hat{\theta})\varphi_{i,j}(x) \right)^{2} \label{eq:pl427}
\end{align}
Defining for all $(x,y,\hat{\theta}) \in \mathbb{R}^{n+1} \times \mathbb{R}^{p}$
\begin{align}
 M(x,y,\hat{\theta}):=& \alpha +\omega+\frac{r+\vert \hat{\theta} \vert}{\omega} \left( \sum_{j=1}^{p}\rho(x,y,\hat{\theta}) \right)^{2}  +\sum_{j=1}^{p} \left( \sqrt{r}+\frac{1}{2}+\frac{1}{2}\hat{\theta}_{j}^{2} \right) \rho(x,y,\hat{\theta}) \nonumber \\
& +\frac{1}{4\epsilon} \sum_{j=1}^{p} \left( \varphi_{n+1,j}(x,y)-\sum_{i=1}^{n}\frac{\partial k}{\partial x_{i}}(x,\hat{\theta})\varphi_{i,j}(x) \right)^{2} \label{eq:pl428}
\end{align}
we obtain inequalities \eqref{eq:inV2}, \eqref{eq:inrt2} from \eqref{eq:pl419}, \eqref{eq:pl428} and definition \eqref{eq:pl4dd}.
\begin{flushright}
$\blacksquare$
\end{flushright}
\vspace{-5mm}
\subsection{Proof of Theorem \ref{thm:2}}
The proof is divided into two parts. In the first part of the proof we prove the existence of the adaptive controller that guarantees inequalities \eqref{eq:th2inV}, \eqref{eq:th2inrt}. In the second part of the proof, we prove all stability properties for the closed-loop system \eqref{eq:csthl26} with \eqref{eq:th2ac0}, \eqref{eq:th2ac}.

\textsl{1st part:} Let arbitrary constants $r \geq 0$, $\alpha, \omega, \epsilon >0$ and $\gamma_{j}>0$, $j=1,\ldots,p$ be given. Applying Lemma \ref{lem4} inductively, it suffices to show that there exist $C^{\infty}$ functions $k:\mathbb{R}^ \times \mathbb{R}^{p} \rightarrow \mathbb{R}$, $w:\mathbb{R} \times \mathbb{R}^{p} \rightarrow \mathbb{R}^{p}$ with $k(0,\hat{\theta})=0$ and $w(0,\hat{\theta})=0$ for all $\hat{\theta} \in \mathbb{R}^{p}$ and a smooth patameterized family of diffeomorphisms of $\mathbb{R}$ $T:\mathbb{R} \times \mathbb{R}^{p} \rightarrow \mathbb{R}$ such that the closed-loop system
\begin{equation}
\dot{x}_{1}=f_{1}(x_{1})+g_{1}(x_{1})x_{2}+\sum_{j=1}^{p} \varphi_{1,j}(x_{1})\theta_{j} \label{eq:pth233}
\end{equation}
with the adaptive feedback law
\begin{align}
\dfrac{d\hat{\theta}}{dt} &= w(x_{1},\hat{\theta}) \label{eq:pth234o} \\
x_{2} &= k(x_{1},\hat{\theta}) \label{eq:pth234}
\end{align}
satisfies the following inequalities for all $(x,\hat{\theta},\theta) \in \mathbb{R} \times \mathbb{R}^{p} \times \mathbb{R}^{p}$:
\begin{align}
\frac{d}{dt} \left( \frac{1}{2}\vert T(x_{1},\hat{\theta}) \vert^{2}+  \sum_{j=1}^{p}\frac{1}{2\gamma_{j}}\left( \hat{\theta}_{j}-\theta_{j} \right)^{2} \right) 
 \leq -\alpha\vert T(x_{1},\hat{\theta}) \vert^{2} \label{eq:pth235} 
\end{align}
\begin{align}
\frac{d}{dt} \left( \frac{1}{2}\vert T(x_{1},\hat{\theta}) \vert^{2} \right) \leq& -\omega \vert T(x_{1},\hat{\theta}) \vert^{2} + 2^{-n}\epsilon \left( \vert \theta \vert^{2}-r \right)^{+} \label{eq:pth236}
\end{align}
Define the map $T:\mathbb{R} \times \mathbb{R}^{p} \rightarrow \mathbb{R}$ by $T(x_{1},\hat{\theta})=x_{1}$ for all $x_{1}\in \mathbb{R}$ and the functions 
\begin{align}
& k(x_{1},\hat{\theta})= -\dfrac{1}{g_{1}(x_{1})}\left( f_{1}(x_{1})+\sum_{j=1}^{p}\varphi_{1,j}(x_{1})\hat{\theta}_{j} + \tilde{M}(x_{1},\hat{\theta})x_{1} \right) \label{eq:pth237} 
\end{align}
\begin{align}
&w(x_{1},\hat{\theta})=x_{1}\left( \gamma_{1}\varphi_{1,1}(x_{1}),\gamma_{2}\varphi_{1,2}(x_{1}),\ldots,\gamma_{p}\varphi_{1,p}(x_{1}) \right)' 
\end{align}
where $\tilde{M}(x_{1},\hat{\theta})$ is determined in what follows. The closed-loop system \eqref{eq:pth233} with the adaptive feedback law \eqref{eq:pth234o}, \eqref{eq:pth234} satisfies the following equations for all $(x,\hat{\theta},\theta) \in \mathbb{R} \times \mathbb{R}^{p} \times \mathbb{R}^{p}$:
\begin{align}
\frac{d}{dt} \left( \frac{1}{2}\vert T(x_{1},\hat{\theta}) \vert^{2}  + \sum_{j=1}^{p}\frac{1}{2\gamma_{j}}\left( \hat{\theta}_{j}-\theta_{j} \right)^{2} \right)  =-\tilde{M}(x_{1},\hat{\theta})x_{1}^{2} \label{eq:pth239}
\end{align}
\vspace{-0mm}
\begin{align}
\frac{d}{dt} \left( \frac{1}{2}\vert T(x_{1},\hat{\theta}) \vert^{2} \right)=& -\tilde{M}(x_{1},\hat{\theta})x_{1}^{2}+ x_{1}\sum_{j=1}^{p}\varphi_{1,j}(x_{1})\theta_{j}-x_{1}\sum_{j=1}^{p}\varphi_{1,j}(x_{1})\hat{\theta}_{j} \label{eq:pth240}
\end{align}
Define the function $\bar{P}:\mathbb{R}^{p} \rightarrow \mathbb{R}^{p}$ by means of \eqref{eq:projdef}. We get from \eqref{eq:pth240} for all $(x_{1},\hat{\theta}) \in \mathbb{R} \times \mathbb{R}^{p}$:
\begin{align}
 \frac{d}{dt} \left( \frac{1}{2}\vert T(x_{1},\hat{\theta}) \vert^{2} \right)=& -\tilde{M}(x_{1},\hat{\theta})x_{1}^{2}+x_{1}\sum_{j=1}^{p}\varphi_{1,j}(x_{1})\bar{P}_{j}(\theta) \nonumber \\
& +x_{1}\sum_{j=1}^{p}\varphi_{1,j}(x_{1})(\theta_{j}-\bar{P}_{j}(\theta))+x_{1}\sum_{j=1}^{p}\varphi_{1,j}(x_{1})\hat{\theta}_{j} \label{eq:pth241}
\end{align}
By virtue of Lemma \ref{lem1}, there exist $C^{\infty}$ maps $\rho_{j}:\mathbb{R} \times \mathbb{R}^{p} \rightarrow \left[ 1,+\infty \right)$, $j=1,\ldots,p$ for which the following inequalities hold for $j=1,\ldots,p$:
\begin{equation}
\vert \varphi_{1,j}(x_{1}) \vert \leq \rho_{j}(x_{1})\vert x_{1} \vert \text{ for all } x_{1} \in \mathbb{R} \label{eq:pth242}
\end{equation}
The inequalities
\begin{align}
x_{1}\varphi_{1,j}(x_{1})&(\theta_{j}-\bar{P}_{j}(\theta)) \leq  2^{-n}\epsilon(\theta_{j}-\bar{P}_{j}(\theta))^{2}+\frac{2^{n-2}}{\epsilon}x_{1}^{2}\varphi_{1,j}^{2}(x_{1})
\end{align}
which hold for all $j=1,\ldots,p$ and the inequalities \eqref{eq:pth242} combined with \eqref{eq:pth241} give for all $(x_{1},\hat{\theta}) \in \mathbb{R} \times \mathbb{R}^{p}$:
\begin{align}
\frac{d}{dt} \left( \frac{1}{2}\vert T(x_{1},\hat{\theta}) \vert^{2} \right) \leq& -\tilde{M}(x_{1},\hat{\theta})x_{1}^{2}+x_{1}^{2}\sum_{j=1}^{p}\rho_{j}(x_{1})\vert \bar{P}_{j}(\theta)\vert +2^{-n}\epsilon \sum_{j=1}^{p}(\theta_{j}-\bar{P}_{j}(\theta))^{2} \nonumber \\ 
& +\frac{2^{n-2}}{\epsilon}\sum_{j=1}^{p}x_{1}^{2}\varphi_{1,j}^{2}(x_{1}) +x_{1}^{2}\sum_{j=1}^{p}\rho_{j}(x_{1})\vert \hat{\theta}_{j}\vert \label{eq:pth243}
\end{align}
Definition \eqref{eq:projdef} implies that $\vert \bar{P}_{j}(\theta) \vert \leq \sqrt{r}$ for $j=1,\ldots,p$ and $\vert \theta-\bar{P}(\theta) \vert^{2} \leq \left( \vert \theta \vert^{2}-r \right)^{+}$. Using the previous inequalities in conjunction with \eqref{eq:pth243} and the inequalities $\vert \hat{\theta}_{j}\vert \leq \frac{1}{2}+\frac{1}{2}\hat{\theta}_{j}^{2}$ for $j=1,\ldots,p$, we obtain for all $(x_{1},\hat{\theta}) \in \mathbb{R} \times \mathbb{R}^{p}$:
\begin{align}
 \frac{d}{dt} \left( \frac{1}{2}\vert T(x_{1},\hat{\theta}) \vert^{2} \right) \leq& 
-\left( \tilde{M}(x_{1},\hat{\theta})-\frac{2^{n-2}}{\epsilon}\sum_{j=1}^{p}\varphi_{1,j}^{2}(x_{1})  -\frac{1}{2}\sum_{j=1}^{p}(1+2\sqrt{r}+\hat{\theta}_{j}^{2})\rho_{j}(x_{1}) \right)x_{1}^{2} \nonumber \\
& +2^{-n}\epsilon \left( \vert \theta \vert^{2}-r \right)^{+} \label{eq:pth244}
\end{align}
Inequalities \eqref{eq:pth235}, \eqref{eq:pth236} are direct consequences of \eqref{eq:pth239}, \eqref{eq:pth244} and definitions $T(x_{1},\hat{\theta}):=x_{1}$ and 
\begin{align}
\tilde{M}(x_{1},\hat{\theta}):=& \alpha + \omega+\frac{2^{n-2}}{\epsilon}\sum_{j=1}^{p}\varphi_{1,j}^{2}(x_{1}) +\frac{1}{2}\sum_{j=1}^{p}(1+2\sqrt{r}+\hat{\theta}_{j}^{2})\rho_{j}(x_{1}) 
\end{align}
for all $(x_{1},\hat{\theta}) \in \mathbb{R} \times \mathbb{R}^{p}$.

\textsl{2nd part}: Lemma \ref{lem3} guarantees that for every $\rho \geq 0$ there exist $K_{\infty}$ functions $a_{\rho,\theta},b_{\rho,\theta}:\mathbb{R}_{+} \rightarrow \mathbb{R}_{+}$ such that 
\begin{equation}
a_{\rho,\theta}(\vert x \vert) \leq \vert T(x,\hat{\theta}) \vert \leq b_{\rho,\theta}(\vert x \vert)  \label{eq:pth245}
\end{equation}
for all $(x,\hat{\theta}) \in \mathbb{R}^{n} \times \mathbb{R}^{p}$ with $\vert \hat{\theta}-\theta \vert \leq \rho$.\\
Define the function
\begin{equation}
W(x,\hat{\theta}):= \left( \left( \vert T(x,\hat{\theta})\vert^{2}-\omega^{-1}\epsilon\left( \vert \theta \vert^{2}-r \right)^{+} \right)^{+}\right)^{2} \label{eq:pth246}
\end{equation}
for all $(x,\hat{\theta},\theta) \in \mathbb{R}^{n} \times \mathbb{R}^{p} \times \mathbb{R}^{p}$. Using \eqref{eq:th2inrt} and definition \eqref{eq:pth246}, we obtain the following differential inequality:
\begin{equation}
\frac{d}{dt}(W(x,\hat{\theta})) \leq -4 \omega W(x,\hat{\theta}), \text{ for all } (x,\hat{\theta}) \in \mathbb{R}^{n} \times \mathbb{R}^{p} \label{eq:pth247} 
\end{equation}
Let $R \geq 0$ be a given (arbitrary) constant and consider the solution of the closed-loop system \eqref{eq:csthl26} with \eqref{eq:th2ac0}, \eqref{eq:th2ac} and (arbitrary) initial condition $(x(0),\hat{\theta}(0)) \in \mathbb{R}^{n} \times \mathbb{R}^{p}$ that satisfies
\begin{equation}
\vert x(0) \vert^{2}+ \vert \hat{\theta}(0)-\theta \vert^{2} \leq R^{2} \label{eq:pth248}
\end{equation}  
The solution is defined for $t \in \left[0,t_{max} \right)$, where $t_{max} \in \left(0,+\infty \right]$ is the maximal existence time of the solution. Using \eqref{eq:defV}, \eqref{eq:pth245} and \eqref{eq:pth248}, we get:
\begin{equation}
V(x(0),\hat{\theta}(0))\leq \frac{1}{2} \left( b_{\rho,\theta}^{2}(R)+\dfrac{R^{2}}{\min_{i=1,\ldots,n}(\gamma_{i})} \right) \label{eq:pth249} 
\end{equation}
The differential inequality \eqref{eq:th2inV} guarantees the following estimate for all $t \in \left[0,t_{max} \right)$:
\begin{equation}
V(x(t),\hat{\theta}(t)) \leq V(x(0),\hat{\theta}(0)) \label{eq:pth250}
\end{equation}
Using definition \eqref{eq:defV} and \eqref{eq:pth249}, \eqref{eq:pth250} we get the following estimate for all $t \in \left[0,t_{max} \right)$:
\begin{equation}
\vert \hat{\theta}(t)-\theta \vert \leq \sqrt{2 \max_{i=1,\ldots,n}(\gamma_{i})V(x(t),\hat{\theta}(t))} \leq \rho \label{eq:pth251}
\end{equation}
where
\begin{equation}
\rho := b_{R,\theta}(R)\sqrt{\max_{i=1,\ldots,n}(\gamma_{i})}+R\sqrt{\dfrac{\max_{i=1,\ldots,n}(\gamma_{i})}{\min_{i=1,\ldots,n}(\gamma_{i})}} \label{eq:pth252}
\end{equation}
Using \eqref{eq:defV}, \eqref{eq:pth245}, \eqref{eq:pth249}, \eqref{eq:pth250} and \eqref{eq:pth251} we get the following estimate for all $t \in \left[0,t_{max} \right)$:
\begin{equation}
\vert x(t) \vert \leq a_{\rho,\theta}^{-1}\left( b_{R,\theta}(R)+\dfrac{R}{\sqrt{\min_{i=1,\ldots,n}(\gamma_{i})}} \right) \label{eq:pth253}
\end{equation}
It follows from \eqref{eq:pth251} and \eqref{eq:pth253} that the solution is indeed bounded on $\left[ 0,t_{max} \right)$. Therefore, $t_{max}=+\infty$ and estimates \eqref{eq:pth251}, \eqref{eq:pth253} hold for all $t \geq 0$. Since both $R \geq 0$ and $(x(0),\hat{\theta}(0)) \in \mathbb{R}^{n} \times \mathbb{R}^{p}$ are arbitrary, we conclude that the closed-loop system \eqref{eq:csthl26} with \eqref{eq:th2ac0}, \eqref{eq:th2ac} is forward complete. Property (vi) is a direct consequence of Theorem 1 in \cite{b5}, definition \eqref{eq:pth246}, inequality \eqref{eq:pth247} and the inequality
\begin{align}
\frac{d}{dt}V(x,\hat{\theta}) & \leq -\alpha \vert T(x,\hat{\theta}) \vert^{2} \nonumber \\ 
& \leq -\alpha \left( \vert T(x,\hat{\theta}) \vert^{2}-\omega^{-1}\epsilon \left( \vert \theta \vert^{2}-r \right)^{+} \right)^{+} \nonumber \\
& \leq -\alpha \sqrt{W(x,\hat{\theta})} 
\end{align} 
which is a consequence of \eqref{eq:th2inV}, and definition \eqref{eq:pth246}.\\
Definition \eqref{eq:pth246} and inequality \eqref{eq:pth247} imply that every solution of the closed-loop system \eqref{eq:csthl26} with \eqref{eq:th2ac0}, \eqref{eq:th2ac} satisfies the following estimate for all $t \geq 0$:
\begin{align}
\left( \vert T(x(t),\hat{\theta}(t)) \vert^{2}-\omega^{-1}\epsilon \left( \vert \theta \vert^{2}-r \right)^{+} \right)^{+}  
 \leq \exp(-2\omega t)\left( \vert T(x(0),\hat{\theta}(0)) \vert^{2}-\omega^{-1}\epsilon \left( \vert \theta \vert^{2}-r \right)^{+} \right)^{+} \label{eq:pth254}
\end{align}
Using \eqref{eq:pth245}, \eqref{eq:pth248}, \eqref{eq:pth254} and \eqref{eq:pth251}, we obtain the following estimates for all $t \geq 0$:
\begin{equation}
\vert x(t) \vert \leq a^{-1}_{\rho,\theta}\left( \exp(-\omega t)b_{R,\theta}(R)+\sqrt{\omega^{-1}\epsilon\left( \vert \theta \vert^{2}-r \right)^{+}} \right) \label{eq:pth255}
\end{equation}
Properties (ii), (iii) and (v) are direct consequences of estimates \eqref{eq:pth251}, \eqref{eq:pth253}, \eqref{eq:pth255}. Property (i) is a consequence of the fact that $f_{i}(0)=0$, $\varphi_{i,j}(0)=0$, $k(0,\hat{\theta})=0$ and $w(0,\hat{\theta})=0$ for all $\hat{\theta} \in \mathbb{R}^{p}$. Finally, property (iv) is a consequence of estimate \eqref{eq:pth253} and LaSalle's theorem (Theorem 4.4 in \cite{b4}).
\begin{flushright}
$\blacksquare$
\end{flushright}

\section{Conclusions}
In this paper we considered nonlinear adaptive control systems satisfying a matching condition. We designed controllers with nonlinear damping terms depending both on the state and the parameter estimate and derived KL-estimates of the performance of the closed-loop system that guarantees Global Uniform Output Asymptotic Stability for a class of Output functions in the absence of PE. We compare the proposed controller with the existing one in the literature. The new design enhances the robustness and the convergence speed of the closed-loop system but presents a higher overshoot. Further study would be at the direction of improving the performance of the proposed controller in the parameter estimation convergence at its true value. To this end, employing delayed measurements of the state would prove to be sufficient.

This paper should be properly situated within the literature. The paper, and its contribution, should be viewed in the context of of the papers \cite{bKr2} and \cite{bKr}, which dealt, for the first time, with, respectively, \emph{transient performance} and \emph{asymptotic performance} in adaptive control. The present paper also deals with improving transient performance, for the plant/regulated state. The difference from \cite{bKr2} is that \cite{bKr2} deals with ``overshoot performance’’ (mainly in the $L_{\infty}$, but also in the $L_{2}$ temporal norms, without guarantees on the rate of decay), whereas the present paper deals with ensuring a ``decay rate performance’’, in the sense of a KL estimate. Both \cite{bKr2} and the present paper improve the transient performance by strengthening the non-adaptive controller, without significant alterations in the update law. The difference is that \cite{bKr2} merely employs a robust linear controller to ensure overshoot performance, whereas the present paper goes a step further with nonlinear damping in both the plant and estimator states to improve the \emph{decay rate performance}.


\begin{thebibliography}{00}


\bibitem{b6} A. Bacciotti and L. Rosier, \emph{Liapunov Functions and Stability in Control Theory}, Springer, 2001.

\bibitem{b16} C. P. Bechlioulis and G. A. Rovithakis, "Robust Adaptive Control of Feedback Linearizable MIMO Nonlinear Systems With Prescribed Performance",  \textit{IEEE Transactions on Automatic Control}, 53, 2008, 2090-2099.

\bibitem{b31}T. Berger, A. Ilchmann and E.P. Ryan,  "Funnel Control of Nonlinear Systems", \textit{Mathematics of Control Signals and Systems}, 33, 2021, 151–194.

\bibitem{b36}M. Chen, X. Liu,  and H. Wang, "Adaptive Robust Fault-Tolerant Control for Nonlinear Systems With Prescribed Performance", \textit{Nonlinear Dynamics} 81,  2015, 1727–1739.

\bibitem{b2} W.-K. Chen, \emph{Linear Networks and Systems.} Belmont, CA, USA: Wadsworth, 1993.

\bibitem{b33}C. Gao, X.-P. Liu, H.-Q. Wang, N.-N. Zhao and L.-B. Wu, "Adaptive Neural Funnel Control for a Class of Pure-Feedback Nonlinear Systems With Event-Trigger Strategy", \textit{International Journal of Systems Science}, 51, 2020, 2307-2325.

\bibitem{b50f}C. Hua, P. Ning and K. Li, "Adaptive Prescribed-Time Control for a Class of Uncertain Nonlinear Systems", \textit{IEEE Transactions on Automatic Control}, 67, 2022, 6159-6166. 

\bibitem{b15} F. Jia, X. Yan, X. Wang, Junwei Lu and Y. Li, "Robust Adaptive Prescribed Performance Dynamic Surface Control for Uncertain Nonlinear Pure-Feedback Systems", \textit{Journal of the Franklin Institute}, 357, 2020, 2752-2772.

\bibitem{b3} I. Karafyllis and J.-P. Jiang, \emph{Stability and Stabilization of Nonlinear Systems}, Springer-Verlag, London (Series: Communications and Control Engineering), 2011.

\bibitem{b23}I. Karafyllis and M. Krstic, "Adaptive Certainty-Equivalence Control With Regulation-Triggered Finite-Time Least-Squares Identification",  \textit{IEEE Transactions on Automatic Control}, 63, 2018, 3261-3275.

\bibitem{b22}I. Karafyllis, M. Kontorinaki and M. Krstic, "Adaptive Control by Regulation-Triggered Batch Least Squares", \textit{IEEE Transactions on Automatic Control}, 65, 2020, 2842-2855.

\bibitem{b5} I. Karafyllis and A. Chaillet, "Lyapunov Conditions for Uniform Asymptotic Output Stability and a Relaxation of Barbălat’s Lemma", \textit{Automatica}, 132, 2021, 109792.

\bibitem{b4} H. K. Khalil, \emph{Nonlinear Systems}, 2\textsuperscript{nd} Edition, Prentice-Hall, 1996.

\bibitem{b51f}P. Krishnamurthy, F. Khorrami and M. Krstic, "Adaptive Output-Feedback Stabilization in Prescribed Time for Nonlinear Systems With Unknown Parameters Coupled With Unmeasured States", \textit{International Journal of Adaptive Control and Signal Processing}, 35, 2021, 184-202. 

\bibitem{bKr2} M. Krstic, P. V. Kokotović and I. Kanellakopoulos, "Transient-Performance Improvement With a new Class of Adaptive Controllers", \textit{Systems and Control Letters}, 21, 1993, 451-461.

\bibitem{b1} M. Krstic, I. Kanellakopoulos, and P. Kokotovic, \emph{Nonlinear and Adaptive Control Design}, Wiley, 1995.

\bibitem{bKr} M. Krstic, "Invariant Manifolds and Asymptotic Properties of Adaptive Nonlinear Stabilizers", \textit{IEEE Transactions on Automatic Control}, 41, 1996, 817-829.

\bibitem{b17} G. Li, Y. Liu, Y. Li and X. Bu, "Adaptive Back-Stepping Control of High-Order Uncertain Nonlinear Systems That a Funnel Control Scheme With Uncertain Dynamics", in the \textit{Proceedings of the 2020 International Conference on Electrical Engineering and Control Technologies (CEECT)}, 2020, 1-8.

\bibitem{b37}Y. Liu, X. Liu and Y. Jing, "Adaptive Neural Networks Finite-Time Tracking Control for non-Strict Feedback Systems via Prescribed Performance", \textit{Information Sciences}, 468, 2018, 29-46.

\bibitem{b28}X. Liu, R. Ortega, H. Su and J. Chu, "Immersion and Invariance Adaptive Control of Nonlinearly Parameterized Nonlinear Systems", \textit{IEEE Transactions on Automatic Control}, 55, 2010, 2209-2214.

\bibitem{b30}X. Liu, H. Wang, C. Gao and M. Chen, "Adaptive Fuzzy Funnel Control for a Class of Strict Feedback Nonlinear Systems", \textit{Neurocomputing}, 241, 2017, 71-80.

\bibitem{b19} A. Loría, R. Kelly and A. R. Teel, "Uniform Parametric Convergence in the Adaptive Control of Mechanical Systems", \textit{European Journal of Control}, 11, 2005, 87-100.

\bibitem{b14} F. Mazenc, M. de Queiroz and M. Malisoff, "Uniform Global Asymptotic Stability of a Class of Adaptively Controlled Nonlinear Systems",  \textit{IEEE Transactions on Automatic Control}, 54, 2009, 1152-1158.

\bibitem{b13} F. Mazenc, M. Malisoff and M. de Queiroz, "Uniform Global Asymptotic Stability of Adaptive Cascaded Nonlinear Systems With Unknown High-Frequency Gains", \textit{Nonlinear Analysis: Theory, Methods and Applications}, 74, 2011, 1132-1145.

\bibitem{b10} J. Orlowski, A. Chaillet and M. Sigalotti, "Counterexample to a Lyapunov Condition for Uniform Asymptotic Partial Stability", \textit{IEEE Control Systems Letters}, 4, 2020, 397-401.

\bibitem{b35}J. Qiu, K. Sun, T. Wang and H. Gao, "Observer-Based Fuzzy Adaptive Event-Triggered Control for Pure-Feedback Nonlinear Systems With Prescribed Performance", \textit{IEEE Transactions on Fuzzy Systems}, 27, 2019, 2152-2162.

\bibitem{b40}J. G. Rueda-Escobedo and J. A. Moreno, "Strong Lyapunov Functions for two Classical Problems in Adaptive Control", \textit{Automatica}, 124, 2021, 109250.

\bibitem{b20} Y. Song, K. Zhao and M. Krstic, "Adaptive Backstepping With Exponential Regulation in the Absence of Persistent Excitation", in the \textit{Proceedings of the 2016 American Control Conference (ACC)} 2016, 739-744.

\bibitem{b21} Y. Song, K. Zhao and M. Krstic, "Adaptive Control With Exponential Regulation in the Absence of Persistent Excitation", \textit{IEEE Transactions on Automatic Control}, 62, 2017, 2589-2596.

\bibitem{b9} E. D. Sontag, \emph{Mathematical Control Theory: Deterministic Finite Dimensional Systems}, 2\textsuperscript{nd} Edition, Springer, 1998.

\bibitem{b7} E. D. Sontag and Y. Wang, "Notions of Input to Output Stability", \textit{Systems and Control Letters}, 38, 1999, 235-248.

\bibitem{b25}H. Taghavifar, Y. Qin and C. Hu, "Adaptive Immersion and Invariance Induced Optimal Robust Control of Unmanned Surface Vessels With Structured/Unstructured Uncertainties", \textit{Ocean Engineering}, 239, 2021, 109792.

\bibitem{b12} A. R. Teel and L. Praly, "A Smooth Lyapunov Function from a Class KL Estimate Involving two Positive Semidefinite Functions", \textit{ESAIM Control, Optimization and Calculus of Variations}, 2002, 313-367.

\bibitem{b11} V. I. Vorotnikov, \emph{Partial Stability and Control}, Birkhauser, Boston, 1998.

\bibitem{b27}L. Wang and C. M. Kellett, "Adaptive Tracking Control via Immersion and Invariance: An (i)ISS Perspective," in the \textit{Proceedings of the 2019 IEEE 58th Conference on Decision and Control (CDC)}, 2019, 7019-7024.

\bibitem{b32}H. Wang, Y. Zou, P. X. Liu and X. Liu, "Robust Fuzzy Adaptive Funnel Control of Nonlinear Systems With Dynamic Uncertainties", \textit{Neurocomputing}, 314, 2018, 299-309.

\bibitem{b52f}W. Yang, G. Cui, Z. Li and C. Tao, "Adaptive Practical Fixed-Time Control for a Class of Nonlinear Systems with Input Saturation", in the \textit{Proceedings of the 39th Chinese Control Conference (CCC)}, 2020, 440-445.  

\bibitem{b18} A. Yarza, V. Santibanez and J. Moreno-Valenzuela, "An Adaptive Output Feedback Motion Tracking Controller for Robot Manipulators: Uniform Global Asymptotic Stability and Experimentation", \textit{International Journal of Applied Mathematics and Computer Science}, 23, 2013, 599-611. 

\bibitem{b53f}J. Yu, A. Stancu, Z. Wu and Y. Wu, "Adaptive Finite/Fixed-Time Stabilizing Control for Nonlinear Systems With Parametric Uncertainty", \textit{International Journal of Robust and Nonlinear Control}, 2022, 1-18.

\bibitem{b34}K. Zhao, Y. Song, C. L. P. Chen and L. Chen, "Adaptive Asymptotic Tracking With Global Performance for Nonlinear Systems With Unknown Control Directions", \textit{IEEE Transactions on Automatic Control}, 67, 2022, 1566-1573.

\bibitem{b54f}K. Zimenko, D. Efimov and A. Polyakov, "Adaptive Finite-Time and Fixed-Time Control Design Using Output Stability Conditions", \textit{International Journal of Robust and Nonlinear Control}, 32, 2022, 6361-6378.

\end{thebibliography}
\end{document}